\documentclass[11pt]{article}

\usepackage{amsmath}
\usepackage{amsfonts}
\usepackage{amsthm}
\usepackage{amssymb}
\usepackage{tikz}
\usepackage{soul}
\usepackage{mathtools}
\usepackage{graphicx}
\usepackage{subfigure}
\usepackage{authblk}
\usepackage{url}

\providecommand{\keywords}[1]{\textbf{\textit{Keywords:}} #1}

\theoremstyle{plain}  
\newtheorem{thm}{Theorem}[section]
\newtheorem{lem}[thm]{Lemma}

\newtheorem{ex}[thm]{Example}

\newtheorem{rem}[thm]{Remark}

\newcommand{\by}{{\mathbf{y}}}

\newcommand{\cA}{{\mathcal A}}

\newcommand{\cF}{{\mathcal F}}

\newcommand{\cM}{{\mathcal M}}

\newcommand{\cT}{{\mathcal T}}

\newcommand{\cX}{{\mathcal X}}

\newcommand{\sG}{{\mathsf G}}

\newcommand{\sJ}{{\mathsf J}}

\newcommand{\sM}{{\mathsf M}}

\newcommand{\F}{{\mathbb{F}}}

\newcommand{\N}{{\mathbb{N}}}
\newcommand{\R}{{\mathbb{R}}}
\newcommand{\Z}{{\mathbb{Z}}}
\newcommand{\PP}{{\mathbb{P}}}

\newcommand{\setof}[1]{\left\{ {#1}\right\}}
\newcommand{\Int}{\mathop{\mathrm{int}}\nolimits}

\newcommand{\length}{\mathop{\mathrm{length}}\nolimits}

\newcommand{\mvmap}{\rightrightarrows}

\newcommand{\pred}[1]{\overleftarrow #1}
\newcommand{\down}{\mathop{{\downarrow}}\nolimits}

\newcommand{\tp}{\text{top}}

    \newfont{\mbf}{msbm10 scaled 1100}

\newcommand{\supp}[1]{\left| {#1}\right|}

\newcommand{\diam}{\mathop{\mathrm{diam}}\nolimits}

\newcommand{\sAtt}{{\mathsf{ Att}}}

\newcommand{\Inv}{\mathop{\mathrm{Inv}}\nolimits}
\newcommand{\Con}{\text{Con}}

\title{Identifying Nonlinear Dynamics with High Confidence from Sparse Data}

\author[,a]{Bogdan Batko\thanks{bogdan.batko@ii.uj.edu.pl}}
\author[,b,c]{Marcio Gameiro\thanks{gameiro@math.rutgers.edu}}
\author[,d]{Ying Hung\thanks{yhung@stat.rutgers.edu}}
\author[,e]{William Kalies\thanks{william.kalies@utoledo.edu}}
\author[,b]{Konstantin Mischaikow\thanks{mischaik@math.rutgers.edu}}
\author[,f]{Ewerton Vieira\thanks{ewerton.vieira@dimacs.rutgers.edu}}

\affil[a]{Division of Computational Mathematics, Faculty of Mathematics and Computer Science, Jagiellonian University, ul. St. Lojasiewicza 6, 30-348 Krak\'{o}w, Poland}
\affil[b]{Department of Mathematics, Rutgers, The State University of New Jersey, Piscataway, NJ, 08854, USA}
\affil[c]{Instituto de Ci\^{e}ncias Matem\'{a}ticas e de Computa\c{c}\~{a}o, Universidade de S\~{a}o Paulo, S\~{a}o Carlos, S\~{a}o Paulo, Brazil}
\affil[d]{Department of Statistics, Rutgers, The State University of New Jersey, Piscataway, NJ, 08854, USA}
\affil[e]{Department of Mathematics and Statistics, University of Toledo, Toledo, OH, 43606, USA}
\affil[f]{DIMACS, Rutgers, The State University of New Jersey, Piscataway, NJ, 08854, USA}

\begin{document}

\maketitle

\begin{abstract}
We introduce a novel procedure that, given sparse data generated from a stationary deterministic nonlinear dynamical system, can characterize specific local and/or global dynamic behavior with rigorous probability guarantees. More precisely, the sparse data is used to construct a statistical surrogate model based on a Gaussian process (GP).  The dynamics of the surrogate model is interrogated using combinatorial methods and characterized using algebraic topological invariants (Conley index). The GP predictive distribution provides a lower bound on the confidence that these topological invariants, and hence the characterized dynamics, apply to the unknown dynamical system (assumed to be a sample path of the GP). The focus of this paper is on explaining the ideas, thus we restrict our examples to one-dimensional systems and show how to capture the existence of fixed points, periodic orbits, connecting orbits, bistability, and chaotic dynamics.
\end{abstract}


\bigskip

\section{Introduction}
\label{sec:intro}

We propose a novel framework, combining topological dynamics and statistical surrogate modeling with uncertainty quantification, through which it is possible to characterize local and global dynamics from data with probability guarantees. Given a data set $\cT = \setof{(x_n,y_n)\in\R^d \times \R^d \mid n=1,\ldots, N}$, where it is assumed that the data is generated by a continuous dynamical system on a compact set $X\subset \R^d$, we give rigorous bounds on the probability that the characterization of the dynamics is correct.

To put our results in context, recall that given a continuous function $f\colon X\to X$ the traditional focus of dynamical systems has been on understanding the structure of \emph{invariant sets}, i.e., subsets $S\subset X$ such that $f(S)=S$, for which fixed points and periodic orbits are simple examples.
On a global level, this is equivalent to understanding the conjugacy classes of $f$, i.e., the set of $g\colon Y\to Y$ such that $h\circ f = g\circ h$ for some homeomorphism $h\colon X\to Y$.
This is impossible in general \cite{foreman:rudolph:weiss}.
Even in more restrictive settings, correctly capturing the invariant sets may require correctly identifying the nonlinearity to an extremely high order of precision, which is often impossible from a given finite data set.
The logistic map and the associated cascade of period doublings is an archetypal example.

For this reason our approach focuses on coarsely characterizing dynamics rather than identifying the underlying nonlinearity.
Our characterization is done via the Conley index, an algebraic topological invariant, from which one can induce the existence of invariant sets and dynamic structure of invariant sets, e.g., existence of fixed points, periodic orbits, heteroclinic orbits, and chaotic dynamics \cite{mischaikow:mrozek}.

\begin{figure}[!htbp]
\centering
\includegraphics[width=0.9\linewidth]{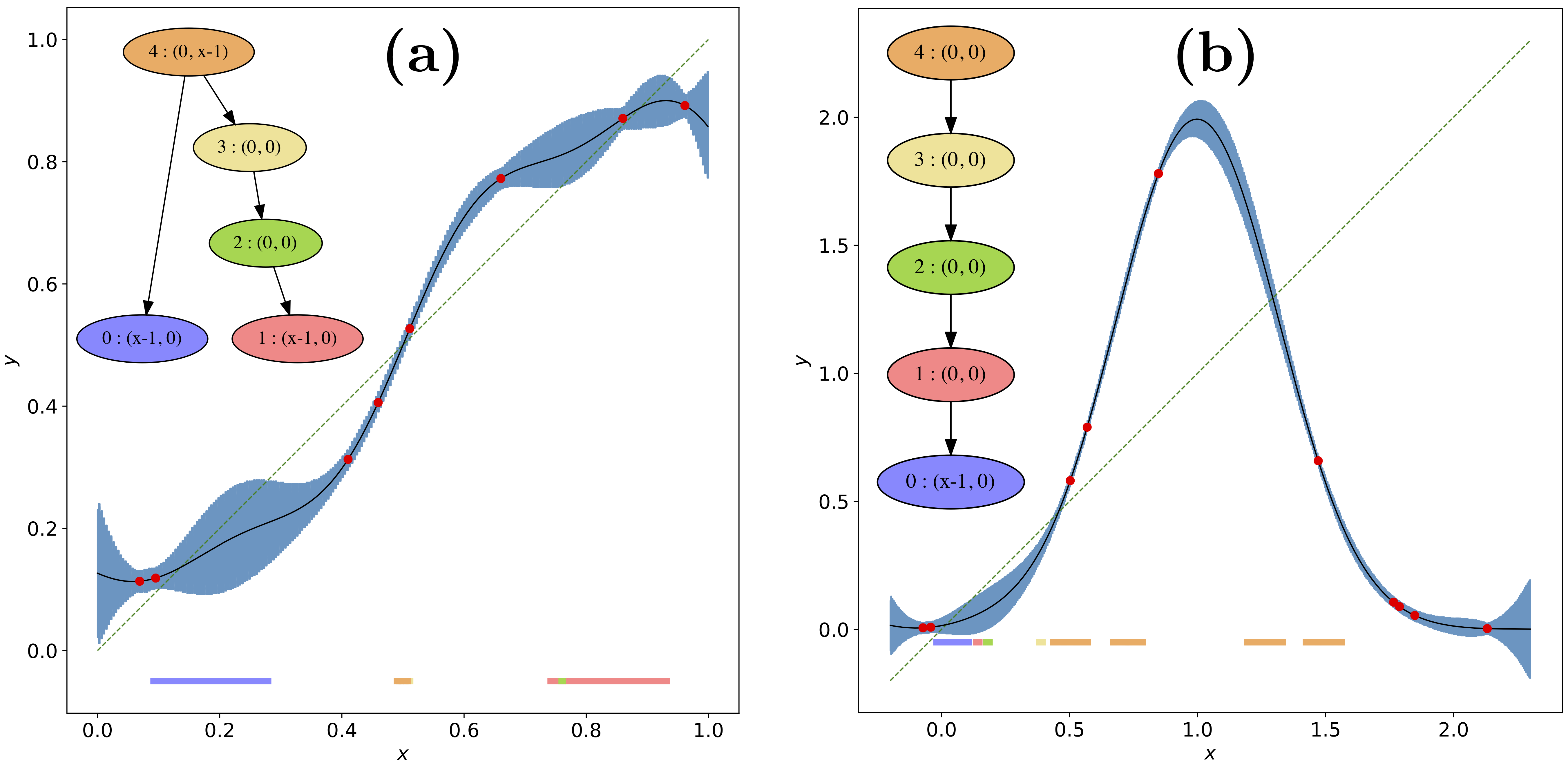}
\caption{In both figures elements of $\cT$ are indicated in red, the mean function $\mu$ is shown in black, and $G$ is shown in blue. The region $G$ is composed of squares of width $2^{-10} \length(X)$. The Morse graphs are indicated at the top left, and the corresponding (color coded) regions  of phase space $\bar{\nu}(\cdot)$ are indicated at the bottom of each figure. {\bf (a)} The Morse graph and the Conley indices of the invariant sets in $\bar{\nu}(\cM(0))$ and $\bar{\nu}(\cM(1))$ -- the blue and red regions on the $x$-axis corresponding to the Morse nodes $0$ (blue) and $1$ (red) -- indicate that bistability is exhibited with $95\%$ confidence. {\bf (b)} The invariant set in $\bar{\nu}(\cM(4))$ -- the orange region on the $x$-axis corresponding to Morse node $4$ -- exhibits chaotic dynamics with positive topological entropy with $95\%$ confidence.}
\label{fig:2examples}
\end{figure}

In essence, the strategy that we propose is straightforward.
It involves a fundamental assumption and three steps that are encapsulated in Fig.~\ref{fig:2examples}.

\begin{description}
\item[A.] Assume the observed data is $\cT = \setof{ (x_n,y_n) \mid y_n = f(x_n) + \epsilon}$, where $f$ is an unknown continuous function. Assume also that there is a Gaussian process (GP) with a prespecified semipositive kernel $k(\cdot,\cdot;\theta)$, where $\theta$ is a vector of unknown parameters associated with the kernel and $\epsilon$ arises from random Gaussian noise, such that $f$ is a realization of this Gaussian process for some value of $\theta$.
\end{description}

\begin{description}
\item[Step 1.] 
Given the data set $\cT$, estimate the unknown parameters and construct a GP surrogate model (see Section~\ref{sec:GP}).
\end{description}

\begin{description}
\item[Step 2.] 
Choose a finite cell complex $\cX$ \cite{lefschetz} whose geometric realization as a regular CW-complex \cite{hatcher} is $X$. Construct a closed set $G \subset X \times X$ with the following property: $G$ is the geometric realization of products of cells from $\cX$ and each fiber $G_x := G\cap \left(\{x\} \times X\right)$, $x\in X$, is nonempty and contractible. Use the combinatorial representation of $G$ to identify potential dynamics and compute their associated Conley indexes (see Section~\ref{sec:conley}).
\end{description}

The set $G$ above represents a (coarse) combinatorial representation of the dynamics as follows: given a pair of cells $\xi \times \xi' \in \cX \times \cX$ whose geometric realization is contained in $G$, we say that $\xi$ maps to $\xi'$ under the combinatorial dynamics (see Sections~\ref{sec:conley} and \ref{sec:surrogate}).

Given a GP $g$ we denote its graph by $\sG(g) := \setof{(x, g(x)) \mid x \in X}$. The Gaussian predictive distribution determines $\PP(\sG(g)\subset G)$.
It is worth emphasizing that $g$ is the Gaussian process (a random variable) and not a realization of the Gaussian process. A \emph{sample path} (or a \emph{realization}) of the GP $g$ is a function $h \colon X \to \R^d$ that is obtained as a realization of $g$ as a random variable (a value of the random variable $g$). Our goal is to compute dynamics which is valid for all sample paths whose graphs are contained in $G$, that is, for all functions in the set
\[
\mathcal{H} := \{ h \colon X \to \R^d \mid h ~\text{is a sample path of}~ g ~\text{and}~ \sG(h)\subset G \}.
\]
We denote $\PP(\mathcal{H})$ by $\PP(\sG(g)\subset G)$. Since the dynamics is computed using the combinatorial representation of $G$, and the Conley index only depends on $G$, the dynamics computed is valid for all functions in $\mathcal{H}$ and $\PP(\sG(g)\subset G)$ provides a lower bound on the probability that the dynamics identified in {\bf Step 2} occurs for the GP $g$. From the discussion in Section~\ref{sec:conley}, in general fibers $G_x$ with smaller diameters lead to greater potential to identify dynamics.

For many applications, the focus is on particular dynamics and/or specific lower bounds on the confidence of the occurrence of the dynamics. Thus, we introduce a third step.

\begin{description}
\item[Step 3.]  Modify $G \subset X\times X$ to both preserve the dynamics of interest and maximize $\PP(\sG(g)\subset G)$.
\end{description}

In this paper we construct $G$ as described in Section~\ref{sec:surrogate}. In this case, the probability $\PP(\sG(g)\subset G)$ provides a \emph{confidence level} that the computed dynamics is valid.

Fig.~\ref{fig:2examples} is meant to provide geometric intuition of {\bf Steps 1} -  {\bf 3}.
In particular, in Fig.~\ref{fig:2examples}~{\bf (a)}, for any sample path $h$ whose graph lies in the blue region we can conclude that the global dynamics generated by $h$ exhibits bistability as well as the existence of at least three fixed points.
For Fig.~\ref{fig:2examples}~{\bf (b)} we can conclude the existence of chaotic dynamics.
In both cases because of the application of {\bf Step 3} we can conclude that the above mentioned dynamics occurs with a confidence of at least $95\%$.

There are three natural questions concerning convergence that arise from the success claimed in  Fig.~\ref{fig:2examples}.
Recall that $f$ is the unknown continuous function that  is assumed to be a realization of the GP and to have generated the data.
The first question is what dynamics of a given function can be identified via the approximation methods (briefly described in Section~\ref{sec:conley}) of {\bf Step 2}?
A precise answer (see \cite[Theorem 1.3]{kalies:mischaikow:vandervorst:15}) goes beyond the scope of this paper.
An imprecise answer is that for many applications most invariant sets of practical interest are identifiable.
The second and third questions are intertwined and address the level of confidence to which our claims on the dynamics can be accepted. That is, how large can we make $\PP(\sG(g)\subset G)$ in {\bf Step 3}? and, to what level of confidence can we approximate $f$ from data?
Theorem~\ref{thm:convergence} of this paper indicates that the confidence level for both questions can be made arbitrarily large simultaneously assuming that $\cT$ contains sufficiently many data points and that the diameter of the elements of geometric realization of $\cX$ are sufficiently small.

As indicated above, the novelty of our approach arises from the combination of GP surrogate modeling and Conley theory.
While each of these topics are well developed, there does not seem to be much overlap of the associated research communities.
With this in mind we provide a minimal description of surrogate modeling by GP (Section~\ref{sec:GP}) and combinatorial Conley theory (Section~\ref{sec:conley}) in the context of maps on $\R^d$.
However, the examples are given using maps on $\R$, as it allows us to demonstrate the results using simple figures. 

\section{Surrogate Modeling by Gaussian Processes}
\label{sec:GP}

A GP model ({\bf Step 1} in our process), also called kriging in geostatistics, is a widely used surrogate model because of its flexibility, nonlinearity, and the capability of uncertainty quantification through the predictive distribution \cite{santner2003design,Gramacy2020}. 

Recall that our data $\cT$ is generated by the unknown realization $f$ of the GP in assumption {\bf A} with random Gaussian noise so that $y_n = f(x_n)+\epsilon$. Let $f_\ell$ denote the $\ell$-th component of $f$, for $\ell=1,\ldots,d$. Then
\begin{equation}
\label{gp}
f_\ell(x)\sim GP(\beta_\ell, \sigma_\ell^2 k(x,x';\theta_\ell)), 
\end{equation}
where $\beta_\ell$ and $\sigma_\ell^2$ are the unknown mean and variance, and the correlation is defined by the kernel $k(x,x';\theta_\ell) = Corr(f_\ell(x), f_\ell(x');\theta_\ell)$ with $k(x,x';\theta_\ell)=k(x',x;\theta_\ell)$ for $x,x'\in X$. For simplicity of exposition, we assume the data is noise-free and therefore $k(x,x;\theta_\ell)=1$. The results in this paper can be easily extended to noisy data by incorporating nugget effects in the kernel function \cite{Gramacy2020}. There are extensive discussions on correlation functions in the literature \cite{St99}. The mean function $\beta$ can be further extended to include regression terms in the mean function which is known as universal kriging \cite{St99,Cressie1993}.

Based on \eqref{gp}, the maximum likelihood estimators (MLEs) for $\beta_\ell$, $\sigma_\ell^2$, and $\theta_\ell$ can be obtained by
\[
\hat{\beta}_\ell=\frac{{\bf 1}^TK^{-1}(\hat{\theta}_\ell)\by_\ell^T}{{\bf 1}^TK^{-1}(\hat{\theta}_\ell){\bf 1}}, \quad
\hat{\sigma}_\ell^2=\frac{1}{N}(\by_\ell^T-\hat{\beta_\ell})^TK^{-1}(\hat{\theta}_\ell)(\by_\ell^T-\hat{\beta}_\ell),
\]
and
\[
\hat{\theta}_\ell=\arg\min_{\theta} \{N \log(\hat{\sigma}_\ell^2)+\log|K({\theta_\ell})|\}
\]
where  $\by_\ell  :=(f_\ell(x_1),\dots,f_\ell(x_N))$, ${\bf 1}$ is a column of 1's with length $N$, $K({\theta})$ is an $N\times N$ correlation matrix with elements $k(x_i,x_j;{\theta})$ for $1\leq i,j \leq N$, and $|K({\theta})|$ is the determinant of $K({\theta})$.
Other estimation approaches, such as the restricted maximum likelihood (REML) method and estimations by cross validation are also applicable \cite{Cressie1993,santner2003design}. Alternatively, assumption {\bf A} can be regarded as a Bayesian prior on the unknown function $f$ and a fully Bayesian approach can be applied to perform estimation and prediction \cite{santner2003design,Gramacy2020}. In this paper, the parameters are estimated by the MLEs.

In {\bf Step 1},
the prediction for an untried $x\in X$ can be obtained by a $d$-dimensional multivariate normal distribution, $MN(\mu(x),\Sigma(x))$, where $\mu(x) = (\mu_1(x), \ldots, \mu_d(x))$ is the best linear unbiased predictor (BLUP) with
\[
\mu_\ell(x)=E(f_\ell(x)|\cT)=\hat{\beta}_\ell+k(x;\hat{\theta}_\ell)^TK^{-1}(\hat{\theta}_\ell)(\by_\ell^T-\hat{\beta}_\ell),
\]
and the covariance matrix $\Sigma(x)$ has diagonal elements 
\[
Var(f_\ell(x)|\cT)=\hat{\sigma}_\ell^2\left( 1-k(x;\hat{\theta}_\ell)^T K^{-1}(\hat{\theta}_\ell) k(x;\hat{\theta}_\ell)\right),
\]
where $k(x;\hat{\theta})$ is the correlation between the new observation and the existing data, i.e., $k(x;\hat{\theta}) = (k(x,x_1;\hat{\theta}), \ldots, k(x,x_N;\hat{\theta}))$, and $K(\hat{\theta})$ is an $N\times N$ correlation matrix with elements $k(x_i,x_j;\hat{\theta})$ for $1\leq i,j \leq N$. The off-diagonal elements in $\Sigma(x)$ are zeros if the $d$-dimensional outputs are assumed to be independent. By further assuming some correlation structures among outputs through the kernel function $k$, the off-diagonal elements can be estimated by techniques such as co-kriging \cite{furrer}. Note that, using different kernel functions, the correlation structure between the $d$-dimensional outputs can be captured \cite{furrer,LeGratiet2015} and nugget effects \cite{ LeeOwen2015,Gramacy2020} can be included into the kernel function to estimate the sampling error or extrinsic noise associated with the observations.  

Recall that for any random variable $Z\sim MN(m,\Lambda)$, the squared Mahalanobis distance
\[
\rho_{\Lambda} ^2(Z,m):=(Z-m)^T\Lambda^{-1}(Z-m)
\]
has $\chi^2$-distribution with $d$ degrees of freedom \cite{santner2003design}.  

Note that if $g$ is a GP on a parameter space $X$, then the above applies to $g(x)$ at each point $x\in X$. Let $\mu$ and $\Sigma$ denote the predictive mean and covariance functions of $g$, respectively. Accordingly, for any $\delta\in (0,1)$ and any fixed $x\in X$, we have
\[
\PP\left(g(x)\in E_{\Sigma(x)}\left(\mu(x),\chi^2 _d(1-\delta)\right)\right)= 1-\delta,
\]
where
$
E_{\Sigma(x)}(\mu(x),c):=\left\{y\in X\ |\ \rho_{\Sigma(x)} ^2(y,\mu(x))<c\right\}
$
is the \emph{confidence ellipsoid},
$\chi^2 _d(1-\delta)$ stands for a $\chi^2$ quantile of order $1-\delta$ with $d$ degrees of freedom, 
and $1-\delta$ is the \emph{confidence level}.
More generally, if $S \subset X$ is finite and $\delta \in(0,1)$, then there exists a 
function $r\colon S\to (0,\infty)$
such that
\begin{equation}
\label{eq:conf}
\PP\left(g(v)\in E_{\Sigma(v)}(\mu(v),r(v))\ \forall v\in S\right)\geq 1-\delta.
\end{equation}
Observe that the function $r$ is not unique.

\begin{figure}[!htbp]
\centering
\includegraphics[width=0.5\linewidth]{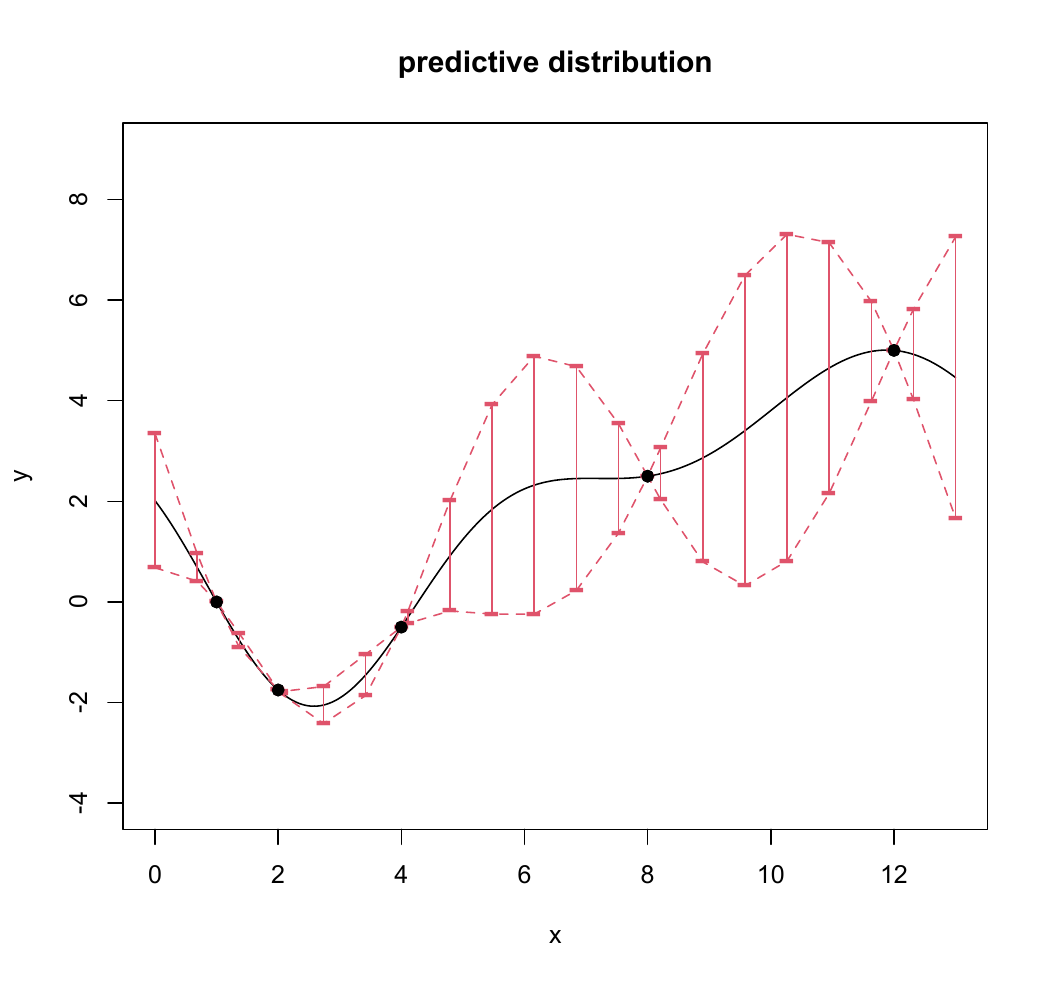}
\caption{Based on five observations, illustrated by solid dots, a one-dimensional GP model is fitted. The black curve is the best linear unbiased predictor and the red bars are the pointwise $95\%$ confidence intervals calculated at $20$ untried points.}
\label{fig:1dGP}
\end{figure}

To illustrate the aforementioned procedure we present a simple one-dimensional ($d=1$) example. Figure~\ref{fig:1dGP} shows a GP model constructed from five noise-free observations indicated by the solid dots. By utilizing a squared exponential kernel we obtain the BLUP $\mu(x)$, depicted as the black curve, which interpolates the observed data. Additionally, the red bars represent the pointwise $95\%$ confidence intervals ($\delta=0.05$) evaluated at $20$ untried points. These confidence intervals quantify the pointwise prediction uncertainty, which decreases to zero when predicting the observed inputs.

In light of {\bf Step 3},
to provide a lower bound on the confidence of our characterization of dynamics, it is reasonable to make use of the pointwise bounds of \eqref{eq:conf} and insist that $G$ satisfies the property that $\left\{y \in X \mid \rho_{\Sigma(x)}^2(y, \mu(v)) \leq r(v)\right\}\subset G_v$ for each $v\in S$.
To obtain appropriate conditions on $G_x$ for $x\in X\setminus S$  we restrict (in this paper) our attention  to kernel functions that are differentiable up to order four, e.g., the squared exponential covariance function or the Mat\'{e}rn kernels with $\nu>2$  \cite{St99}, in which case 
there exists $L_0 > 0$ and constants $a,b>0$ such that for any $L > L_0$ we have
\begin{equation}
\label{eq:L}
\PP\left(\forall x_1,x_2\in X\   \|g(x_1)-g(x_2)\|\leq L\|x_1-x_2\|\right)> \left( 1-ae^{-(\frac{L}{b})^2} \right)^{d^2}
\end{equation}
(cf.\ \cite[Theorem 5]{GhRo06}). 

\section{Combinatorial Conley Theory and the Characterization of Dynamics}
\label{sec:conley}

There are three essential components of combinatorial Conley theory: a finite combinatorial representation of phase space via a cell complex, a combinatorial representation of dynamics via a directed graph, and homological computations.
As described at the end of this section, the combinatorial theory is used to characterize the dynamics generated by continuous functions that are sample paths of the GP.
To instantiate these ideas throughout this section we describe a particularly simple example that leads to bistability.

Recall  \cite{lefschetz} that a \emph{cell complex}  $\cX = (\cX,\leq,\dim,\kappa)$ 
is a finite partially ordered set (poset)  $(\cX,\leq)$, where the partial order $\leq$ indicates the face relation, together with two associated functions \emph{dimension}, $\dim\colon \cX\to \N$, and \emph{incidence}, $\kappa\colon \cX\times \cX\to \F$,
where $\F$ is a principal ideal domain,
subject to the following conditions for all $\xi,\xi',\xi''\in\cX$: (i) $\xi\leq \xi'$ implies $\dim \xi \leq \dim \xi'$, (ii)  $\kappa(\xi',\xi)\neq 0$ implies $\xi'\leq \xi$ and $\dim(\xi) = \dim(\xi')+1$, and $\sum_{\xi'\in \cX} \kappa(\xi,\xi') \kappa(\xi',\xi'')=0$.
A cell complex generates a \emph{chain complex} that we denote by $C_*(\cX;\F)$.
An element $\xi\in \cX$ is called a {\em cell}, and we denote the maximal elements of $\cX$ by $\cX^\tp$.
Given $\xi\in\cX^\tp$, define $\down(\xi) = \setof{\xi'\in \cX\mid \xi'\leq \xi}$.
Note that given any subset $\cA\subset \cX^\tp$, $\down(\cA)$ generates a chain complex $C_*(\cA;\F)$.

\begin{figure}[!htbp]
\centering
\includegraphics[width=1.0\linewidth]{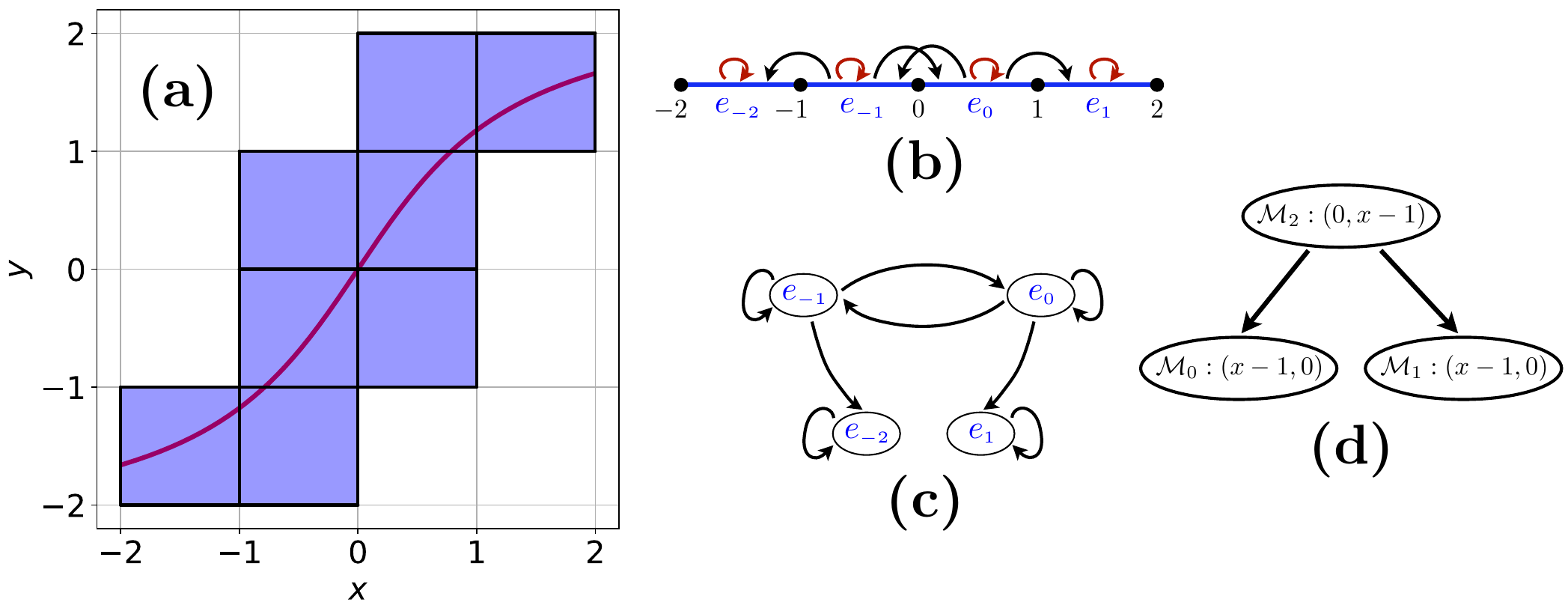}
\caption{Simple example exhibiting bistability. {\bf (a)} The map $f \colon [-2,2] \to [-2,2]$ given by $f(x) = \frac{3}{2}\arctan(x)$ is plotted in red. The domain $X = [-2,2]$ is decomposed into unit intervals $[i, i + 1]$, $i = -2, -1, 0, 1$, and the graph of $f$ is covered by products of these intervals. {\bf (b)} The covering of the graph of $f$ in {\bf (a)} can be represented by a (multi-valued) cell mapping $\cF$ on the set of edges $\setof{e_{-2}, e_{-1}, e_0, e_1}$, where the edge $e_i$ represents the interval $[i, i + 1]$. {\bf (c)} Directed graph representation of the multi-valued map $\cF$. The non-trivial strongly connected components (recurrent components) of this graph give the nodes of the Morse graph $\cM_0 = \setof{e_{-2}}$, $\cM_1 = \setof{e_{1}}$, and $\cM_2 = \setof{e_{-1},e_0}$. {\bf (d)} Morse graph with the Conley index of each node.}
\label{fig:simple_example}
\end{figure}

\begin{ex}
\label{ex:trivial}
Consider $\cX = \setof{v_i\mid i=-2,-1,0,1,2} \cup \setof{e_i\mid i=-2,-1,0,1}$ where we define $\dim(v_i) = 0$ and $\dim(e_i) = 1$, i.e., $v_i$ is a vertex and $e_i$ is an edge (see Figure~\ref{fig:simple_example}). We set $v_i \leq e_j$, if $i=j$ or $i=j+1$. We assume that $\F$ is the field $\Z_2$ and set $\kappa(v_i,e_j) = 1$ if and only if $v_i \leq e_j$. This gives rise to the chain complexes $C_0(\cX;\Z_2)\cong \Z_2^5$ and $C_1(\cX;\Z_2)\cong \Z_2^4$ with a boundary operator $\partial_1\colon C_1(\cX;\Z_2)\to C_0(\cX;\Z_2)$ given by the matrix with entries $\kappa(v_i, e_j)$.
\end{ex}

We represent dynamics using a \emph{combinatorial multivalued map} $\cF\colon \cX^\tp \mvmap \cX^\tp$, i.e., for each $\xi\in\cX^\tp$, $\cF(\xi)\subset \cX^\tp$.
A combinatorial multivalued map is  equivalent to a directed graph with vertices $\cX^\tp$ and edges $\xi\to\xi'$ if $\xi'\in \cF(\xi)$.
To identify the potential recurrent and  gradient-like structure of $\cF$, we make use of the \emph{condensation graph} of $\cF$  obtained by identifying each strongly connected component of $\cF$ to a single vertex \cite{cormen:leiserson:rivest:stein}.
As this is a directed acyclic graph, it can be viewed as a poset that we denote by $\mathsf{SC(\cF)}$.
A \emph{recurrent component} is a strongly connected component that contains at least one edge.
The \emph{Morse graph} of $\cF$, denoted by $\sM(\cF)$, is the subposet of  recurrent components of $\mathsf{SC(\cF)}$.
We typically display the Morse graph via the Haase diagram of $\sM(\cF)$.
Observe that the order relation on $\sM(\cF)$ provides a combinatorial description of the structure of the gradient-like dynamics.

\begin{ex}
\label{ex:trivial2}
Continuing with Example~\ref{ex:trivial}, $\cX^\tp = \setof{e_i\mid i=-2,-1,0,1}$.
Set $\cF\colon \cX^\tp\mvmap \cX^\tp$ to be
\[
\cF(e_{-2})= \setof{e_{-2}},\quad
\cF(e_{-1})= \setof{e_{-2},e_{-1},e_{0}},\quad
\cF(e_{0})= \setof{e_{-1},e_{0},e_{1}},\quad\text{and}\quad
\cF(e_{1})= \setof{e_{1}}.
\]
The strongly connected components are $\cM_0 = \setof{e_{-2}}$, $\cM_1 = \setof{e_{1}}$, and $\cM_2 = \setof{e_{-1},e_0}$. The condensation graph has edges $\cM_2 \to \cM_0$ and $\cM_2 \to \cM_1$. This is an acyclic directed graph and thus can be thought of as a poset with relations $\cM_0 < \cM_2$ and $\cM_1 < \cM_2$. Note that each strongly connected component has at least one edge and therefore each strongly connected component is a recurrent component. Thus the Morse graph $\sM(\cF)$ is the poset with elements $\setof{\cM_0, \cM_1, \cM_2}$. The dynamics interpretation is that one can move from state $\cM_2$ to state $\cM_0$ or to state $\cM_1$, but one cannot move from state $\cM_0$ or state $\cM_1$ to any other state.
\end{ex}

An alternative perspective for characterizing the dynamics associated with $\cF$ is to consider its  \emph{attractors} defined by  $\sAtt(\cF):=\setof{\cA\subset \cX^\tp \mid \cF(\cA) = \cA}$.
The equivalence arises from the fact that, as shown in \cite{kalies:mischaikow:vandervorst:15},  $\sAtt(\cF)$ is a bounded distributive lattice where the partial order is inclusion.
More precisely, if we let $\sJ(\sAtt(\cF))$ denote the set of join irreducible elements of $\sAtt(\cF)$, i.e., those elements of $\sAtt(\cF)$ that have a unique immediate predecessor under inclusion, then there exists a poset isomorphism $\nu\colon \sM(\cF) \to \sJ(\sAtt(\cF))$ \cite{kalies:mischaikow:vandervorst:21}.

\begin{ex}
\label{ex:trivial3}
Continuing with Example~\ref{ex:trivial2}, 
\[
\sAtt(\cF) = \setof{ \emptyset,\setof{e_{-2}},\setof{e_{1}},\setof{e_{-2},e_1},\setof{e_{-2},e_{-1},e_0,e_1}}.
\]
Observe that $\sJ(\sAtt(\cF)) = \setof{ \setof{e_{-2}},\setof{e_{1}},\setof{e_{-2},e_{-1},e_0,e_1}}$     
and using inclusion to define the partial order we obtain a poset that is isomorphic to $\sM(\cF)$.
In this case the poset isomorphism $\nu\colon \sM(\cF) \to \sJ(\sAtt(\cF))$ satisfies $\nu(\cM_0) = \setof{e_{-2}}$, $\nu(\cM_1) = \setof{e_{1}}$, and $\nu(\cM_2) = \setof{e_{-2},e_{-1},e_0,e_1}$. 
\end{ex}

We use the first perspective (associated with posets, e.g., Morse graphs) for efficient computations and to organize the global information, and the second perspective (associate with lattices, e.g., attractors) to identify the homological computations that recover nontrivial information about the structure of the dynamics exhibited by the continuous function.

For the sake of simplicity we define an \emph{index pair}  for $\cF$ to be a pair  $\cA = (\cA_1,\cA_0)$ where $\cA_1,\cA_0\in\sAtt(\cF)$ and $\cA_0\subset\cA_1$. Observe that $\cF(\cA_i)\subset \cA_i$, $i=0,1$. Under rather weak conditions \cite{harker:kokubu:mischaikow:pilarczyk} (we return to this point below) $\cF$ induces a map on homology, i.e.,
\[
\cF_*\colon H_*\left(\down(\cA_1),\down(\cA_0);\F\right)\to H_*(\down(\cA_1),\down(\cA_0);\F).
\]
The \emph{Conley index} of $\cA$, denoted by $\Con_*(\cA;\F)$, is defined to be the shift equivalence class of $\cF_*$ (if $\F$ is a field, then this  is equivalent to the rational canonical form of the linear map $\cF_*$ \cite{bush:cowan:harker:mischaikow}). In particular, we can assign a Conley index to each $\cM\in \sM(\cF)$, by declaring $\Con_*(\cM;\F) \cong \Con_*\left(\down(\nu(\cM)),\down(\pred{\nu(\cM)});\F\right)$ where $\pred{\nu(\cM)}$ is the unique immediate predecessor of $\nu(\cM)$.
Given $\cX$ and $\cF$ there exists software \cite{CMGDB} to compute $\sM(\cF)$,  $\sAtt(\cF)$, and $\Con_*(\cM;\F)$ (for this paper we take $\F \cong\Z_5)$.
This software is based on what are essentially combinatorial algorithms and as a consequence are extremely efficient.

\begin{ex}
\label{ex:trivial4}
Continuing with Example~\ref{ex:trivial3}, the set of index pairs are
\[
(\setof{e_{-2}},\emptyset),\quad (\setof{e_{1}},\emptyset),\quad \text{and}\quad (\setof{e_{-2},e_{-1},e_0,e_1},\setof{e_{-2},e_1}).
\]
The attractors are defined in terms of elements of $\cX^{top}$.
To compute homology we need to work with the chain complexes associated with the attractors, i.e.,
\[
\down(\setof{e_{-2}}) = \setof{ e_{-2},v_{-2},v_{-1}}, \quad \down(\setof{e_{1}}) = \setof{ e_{1},v_{1},v_{2}},\quad\text{and}\quad
 \down(\setof{e_{-2},e_{-1},e_0,e_1}) = \cX.
\]
The relative homology of these index pairs are
\[
H_k(\down(\setof{e_{-2}}), \emptyset; \F) \cong H_k(\down(\setof{e_{1}}), \emptyset; \F) \cong
\begin{cases}
\Z_2, & \text{if $k=0$}\\
0, & \text{otherwise}
\end{cases}
\]
and
\begin{align*}
H_k(\down(\setof{e_{-2},e_{-1},e_0,e_1}), \down(\setof{e_{-2},e_1}); \F) & \cong H_k(\setof{e_{-1}, e_0, v_0}, \emptyset; \F) \\
& \cong H_k(\setof{e_{-1}}, \emptyset; \F) \cong \begin{cases}
\Z_2, & \text{if $k=1$}\\
0, & \text{otherwise}.
\end{cases}    
\end{align*}

The  induced maps on homology are the identity maps and thus the rational canonical forms of $(\cF_0, \cF_1)$ are $(x-1,0)$, $(x-1,0)$, and $(0,x-1)$, respectively.
These are the Conley indices of the elements $\cM_0$, $\cM_1$, and $\cM_2$ of the Morse graph, respectively.
\end{ex}

Before relating the above mentioned combinatorial framework to continuous dynamics we recall the following concepts.
Let $g\colon X\to X$ be a continuous map on a compact space.
Given $N\subset X$, the \emph{maximal invariant set} contained in $N$ is given by 
\[
\Inv(N,g) := \{x\in N \mid \exists\, \sigma\colon \Z\to N\ \text{such that $\sigma(0)=x$ and $\sigma(n+1) = g(\sigma(n))$ for all $n\in\Z$}\}.
\]
A compact set $N\subset X$ is an \emph{attracting block} if $g(N)\subset \Int(N)$ and an \emph{isolating neighborhood} if $\Inv(N,g)\subset \Int(N)$.
It is easily checked that an attracting block is an isolating neighborhood.

As suggested above, the phase space $X$ for the dynamics generated by $g$ is represented by the cell complex $\cX$.
In particular, we assume that $X$ is a regular CW complex \cite{hatcher, kalies:mischaikow:vandervorst:15}, and we use the map $\supp{\cdot}\colon \cX \to X$ to identify how the cell complex $\cX$ realizes the regular CW-complex $X$, i.e., given $\xi\in\cX$, if $\dim(\xi) = n$ then $\supp{\xi}$ represents the corresponding regular closed cell in the $n$-skeleton of $X$ and $\supp{\cX} = X$.
In applications, we start with the space $X$ and choose a decomposition $\cX$.
We define $G=\bigcup_{\xi\in\cX}\supp{\xi}\times \supp{\cF(\xi)}$.

\begin{ex}
\label{ex:trivial5}
Returning to Example~\ref{ex:trivial}, note that $\cX$ represents a decomposition of the interval $[-2,2]\subset \R$ where
$\supp{v_i} = i$ and $\supp{e_i} = [i,i+1]$ for all $i$.
This in turn implies (see Example~\ref{ex:trivial2}) that
\[
G = [-2,-1]\times [-2,-1] \cup [-1,0]\times [-2,1] \cup [0,1]\times [-1,2]\cup [1,2]\times [1,2] \subset [-2,2]^2.
\]
Since each vertical fiber of $G$ is an interval, $\cF$ is acyclic for each $\xi$.
\end{ex}

To relate the combinatorial multivalued map $\cF\colon \cX^{\text{top}}\mvmap \cX^{\text{top}}$ with the continuous function $g\colon X\to X$ we make two assumptions.
First, that $\cF$ is an \emph{outer approximation} of $g$, that is, $g(\supp{\xi})\subset \Int(\supp{\cF(\xi)})$ for all $\xi\in \cX^\tp$.
Second, if we extend $\cF$ to all of $\cX$ by setting $\cF(\xi):= \down\left(\setof{\cF(\xi')\mid \xi'\in \cX^\tp,\ \xi\leq \xi'}\right)$, then $\cF(\xi)$ is  \emph{acyclic}, i.e., $\bar{H}_*(\down(\cF(\xi))=0$, where $\bar{H}$ denotes reduced homology. In this case we say that $\cF$ an \emph{acyclic outer approximation} of $g$.

\begin{ex}
\label{ex:trivial6}
Consider $f \colon [-2,2]\to[-2,2]$ defined by $f(x) = \frac{3}{2}\arctan(x)$.
Returning to Example~\ref{ex:trivial5} we let the reader check that (see Figure~\ref{fig:simple_example})
\[
G(f) = \setof{(x,f(x)\mid -2\leq x\leq 2} \subset \Int(G) \subset [-2,2]^2.
\]
This implies that $\cF$ is an acyclic outer approximation of $f$.
\end{ex}

Under these assumptions,  if $\cM\in \sM(\cF)$, then
\begin{equation}
\label{eq:ConleyIndex}
\Con_*\left(\Inv\left(\supp{\nu(\cM)\setminus \pred{\nu(\cM)}},g \right);\F \right) \sim \Con_*\left(\cM;\F\right)
\end{equation}
where $\Con_*$ on the left denotes the classical homology Conley index for maps \cite{mischaikow:mrozek} and $\Con_*$ on the right the Conley index defined above.
As indicated in the introduction, knowledge of the Conley index provides information about the structure of the dynamics of \[\Inv\left(\supp{\nu(\cM)\setminus \pred{\nu(\cM)}},g \right),\] i.e., the computations outlined in this section provide information about the invariant dynamics contained in $\bar{\nu}(\cM):= \supp{\nu(\cM)\setminus \pred{\nu(\cM)}}$.

\begin{ex}
\label{ex:trivial7}
Combining the information from the previous examples we have derived the following information concerning the global structure of the dynamics generated by the map $f \colon [-2,2]\to[-2,2]$. The fact that the Conley indices of $\cM_i$ are not trivial (Example~\ref{ex:trivial4}) implies that $\Inv(\bar{\nu}(\cM_0), f) = \Inv([-2,-1], f) \neq \emptyset$, $\Inv(\bar{\nu}(\cM_1), f) = \Inv([1,2], f) \neq \emptyset$, and $\Inv(\bar{\nu}(\cM_2), f) = \Inv([-1,1], f) \neq \emptyset$. Furthermore, the fact that the Conley index is the identity map on a one-dimensional vector space implies that each of these invariant sets contains a fixed point. Finally, the dynamics of $f$ exhibits bistability since the poset structure on the Morse graph implies that if $x \in [-2,-1]$, then $f^n(x) \in [-2,-1]$ for all $n \geq 1$, and similarly if $x \in [1,2]$, then $f^n(x) \in [1,2]$ for all $n \geq 1$.
\end{ex}

\begin{rem}
The presentation of  Example~\ref{ex:trivial}-\ref{ex:trivial7}  clearly was chosen to follow (and hopefully enlighten upon) the curt review of combinatorial Conley theory. In practice, e.g., in the examples of Section~\ref{sec:examples}, the order of development is different (see Figure~\ref{fig:simple_example}). The first steps involve the choice of the phase space $X$ and the identification of a surrogate model. For this paper, we restrict our attention to $X\subset \R$ (see \cite{vieira23} for applications of these ideas in the context of robotic control where $X\subset \R^n$ with $n>1$). Furthermore, unlike our choice in Example~\ref{ex:trivial} of $\cX^{\tp}$ consisting of four elements, in the examples of Section~\ref{sec:examples} we choose $\cX^{\tp}$ containing at least $2^{10}$ elements. In this setting an explicit list of the rectangular regions that make up the region $G\subset X\times X$ is meaningless, thus we plot $G$ in blue (see Figure~\ref{fig:2examples}). We also present the $\sM(\cF)$ as a graph where the relative ordering decreases as one goes from top to bottom, e.i., minimal elements are at the bottom, and indicate the Conley index within each node.
\end{rem}

\section{Probabilistic Bounds Using Gaussian Process Surrogates}
\label{sec:surrogate}

Let $X\subset\R^d$ be a compact regular CW complex indexed by a cell complex $\cX$, i.e., if $\xi\in \cX$ and $\dim(\xi)=\ell$, then $\supp{\xi}$ is the closure of the $\ell$-dimensional cell in $X$.
We assume that {\bf Step 1} and {\bf Step 2} have been completed, which implies that we have obtained a predictive mean $\mu \colon X \to \R$, a predictive covariance function $\Sigma$ based on the GP, and identified dynamics via Conley theory.
We adopt the approach that we are only interested in dynamics in $X$ that can be obtained with a given amount of confidence quantified by $0<\delta<1$.

Let $S = \cX^{(0)}$, where $\cX^{(0)}$ denotes the set of vertices of $\cX$.
Choose $r(v)>0$ for $v \in S$ such that
\eqref{eq:conf} is satisfied.
We emphasize that there is considerable freedom in the choice of the individual values of $r(v)$. 
In particular, when $r(v) = \chi_d^2 \left( (1-\delta)^{1 / \# S} \right)$ for each $v \in S$ we say that we are choosing \emph{pointwise equal confidence}.

To define $\cF\colon \cX\mvmap\cX$ we make use of the following notation. Let $\diam{A} = \sup_{x,x'\in A}\|x-x'\|$ for $A\subset X$. Given $\xi\in\cX$, set $\diam{\xi} = \diam\supp{\xi}$, and $\diam{\cX} = \sup_{\xi\in\cX}\diam\xi$. Moreover, given $U\subset \R^d$, let $B_z(U) := \{ x\in \R^d \mid \inf_{u\in U}\|x-u\|\leq z\}$.
Choose $L>0$ sufficiently large (see Lemma~\ref{lem:approx}). If $\xi\in\cX^{\text{top}}$, set
\begin{equation}
\label{eq:Q}
Q(\xi) := \bigcap_{v<\xi,\, v\in S} B_{L\diam{\xi}}\left(E_{\Sigma(v)}(\mu(v),r(v))\right)
\end{equation}
and define $\cF\colon \cX^{\text{top}}\mvmap \cX^{\text{top}}$ by $\cF(\xi) := \{\xi'\in \cX^{\text{top}}\mid \supp{\xi'}\cap Q(\xi)\neq\emptyset \}$. Define
\begin{equation}
\label{eq:defnG}
G = G_{L,r} := \bigcup_{\xi\in\cX^{\text{top}}} \supp{\xi}\times\supp{\cF(\xi)} \subset X \times X,
\end{equation}
and
\begin{equation}
\label{eq:defnGtilde}
\tilde{G} = \tilde{G}_{L,r} := G_{L,r} \cup \left( \bigcup_{\xi \in \cX^{\text{top}}} \supp{\xi} \times Q(\xi) \right) \subset X \times \R^d,
\end{equation}
where when convenient we drop the explicit dependence on $L$ and $r$. Note that $G$ is a cover by cells of the cell decomposition of $X \times X$ of the confidence sets given by \eqref{eq:Q} restricted to $X \times X$, while $\tilde{G}$ includes the portions of these confidence sets that are not in $X \times X$, and $G = \tilde{G} \cap (X \times X)$.

\begin{thm}
\label{thm:convergence}
Let $\cT$ be a data set that satisfies assumption {\bf A} where $\setof{x_n\mid n=1,\ldots, N}$ are chosen i.i.d. from a uniform distribution and assume the kernel $k$ satisfies the conditions for \eqref{eq:L}. Let $\alpha>0$ and $\delta\in(0,1)$. There exist $\varepsilon_0>0$ and $n_0\in\N$ such that the set $\tilde{G}$ given by \eqref{eq:defnGtilde} satisfies
\begin{equation}
\label{eq:diam}
\PP\left(\sup_{x\in X}\diam(\tilde{G}_x)<\alpha \right)>1-\delta
\quad\text{and}\quad
\PP( \sG(g) \subset \tilde{G})>1-\delta
\end{equation}
provided that $N > n_0$, $g$ is a GP constructed as in {\bf Step 1}, and as in {\bf Step 2}, $X\subset\R^d$ is a compact regular CW-complex indexed by a cell complex $\cX$ with $\diam(\cX)<\varepsilon_0$, where the top dimensional cells $\cX^\text{top}$ are $d$-dimensional.
\end{thm}

\begin{rem}
Note that $G \subset X \times X$ is defined in terms of the map $\cF$ and so the computed dynamics is valid for all samples paths whose graphs are in $G$ and $\PP( \sG(g) \subset G)$ gives the confidence level on the dynamics. However we can only estimate $\PP( \sG(g) \subset \tilde{G})$, by Theorem~\ref{thm:convergence}, and hence we adopt the following strategy: If $\tilde{G} \subset X \times X$, and hence $G = \tilde{G}$, then we have the confidence level $\PP( \sG(g) \subset G) = \PP( \sG(g) \subset \tilde{G}) > 1 - \delta$ on the dynamics. If, on the other hand, $\tilde{G} \not\subset X \times X$, then we cannot estimate the confidence level $\PP( \sG(g) \subset G)$ and so we declare failure in identifying the dynamics with confidence level $1-\delta$, since in this case $\PP( \sG(g) \subset G) \leq \PP( \sG(g) \subset \tilde{G})$ and hence the confidence level may be less than $1 - \delta$.

Notice that Theorem~\ref{thm:convergence} indicates that if $f(X) \subset \Int(X)$, then we should have $\tilde{G} \subset X \times X$ as long as we have enough data points and the grid is sufficiently fine. Therefore a failure suggests that it may be necessary to choose a larger domain $X$, more data points, or a smaller confidence level $1 - \delta$.
\end{rem}

Theorem~\ref{thm:convergence} implies that with sufficient data and sufficient computational effort we can obtain the following two fundamental results:
\begin{enumerate}
\item Detailed dynamics can be extracted from $G$ via the Conley theory computations, since the sizes of the fibers of $G$ are bounded above by an arbitrarily chosen $\alpha$ with confidence $1-\delta$;
\item The dynamics identified via $G$ occurs since a given realization $h$ of the GP model is a selector of $G$, that is $\sG(h) \subset G$, with probability $1-\delta$. Furthermore, this dynamics is valid with confidence greater than $1-\delta$.
\end{enumerate}
In the computations in Section~\ref{sec:examples} we fix the data size $N$, and hence we only give the confidence level of the correctness of the dynamics (item 2 above).

To prove Theorem~\ref{thm:convergence} we begin by recalling and establishing the necessary notation, and then proving a series of lemmas that are used in the proof.

For the remainder of this section we assume that $X\subset \R^d$ is a compact set that is the regular CW complex realization of a cell complex $\cX$ where the top dimensional cells  $\cX^\text{top}$ are $d$-dimensional.

We assume that $\cT = \setof{(x_n,y_n)\mid n=1,\ldots, N}$ is a data set that satisfies assumption {\bf A}, the kernel $k$ satisfies the conditions for equation \eqref{eq:L}, $g$ is a GP constructed as in {\bf Step 1}, and $\tilde{G} = \tilde{G}_{L,r} \subset X \times \R$ is constructed according to equation \eqref{eq:defnGtilde}, where the conditions on $L$ and $r$ are described in what follows.

Recall that $\sG(g)$ denotes the graph of $g$. We use the following lemma to quantify the confidence that the graph of the GP $g$ lies in $\tilde{G}$.

For Theorem~\ref{thm:convergence} and Lemma~\ref{lem:Gbound} we need the assumption that $S = \cX^{(0)}$, however for the next two Lemmas we have more flexibility on the choice of $S$ as stated. Note that Lemma~\ref{lem:approx} gives the confidence level for our computations.

\begin{lem}
\label{lem:approx}
Fix $\delta\in(0,1)$ and let $S \subset \cX^{(0)}$ be such that for each $\xi\in \cX^\text{top}$, there exists $v \in S$ such that $v < \xi$. Then, there exists $r\colon S\to(0,\infty)$ that satisfies equation \eqref{eq:conf} and $\tilde{L} > 0$ such that   $L > \tilde{L}$ implies that
\begin{equation}\label{eq:incl}
\PP\left(\sG(g)\subset \tilde{G}_{L,r}\right)>1-\delta.
\end{equation}
\end{lem}

\begin{proof}
Let $L_0, a, b > 0$ be as in equation \eqref{eq:L} and let $\delta = \delta_1 + \delta_2$ with $\delta_1, \delta_2 > 0$. Choose $\tilde{L} \geq L_0$ such that $1 - a e^{-(L/b)^2} \geq 1 -\delta_1$ for $L > \tilde{L}$. By equation \eqref{eq:L} with probability greater than $1 - a e^{-(L/b)^2}$ we have
\begin{equation}
\label{eq:L_recall}
\|g(x)-g(y)\|\leq L\|x-y\|\mbox{ for all }x,y\in X.
\end{equation}
Let $\cX$ be a cell complex decomposition of $X$. According to equation \eqref{eq:conf} we can pick $r \colon S \to (0,\infty)$ such that with probability greater than $1 - \delta_2$ we have
\begin{equation}
\label{eq:beta}
g(v)\in E_{\Sigma(v)}(\mu(v),r(v))\text{ for all }v\in S.
\end{equation}
It follows that $g$ satisfies \eqref{eq:L_recall} and \eqref{eq:beta} with probability greater than $1 - \delta$.

Now it suffices to show that if $g$ satisfies \eqref{eq:L_recall} and \eqref{eq:beta}, then $\sG(g)\subset \tilde{G}_{L,r}$. If $x\in X$, then $x\in\supp{\xi}$ for some $\xi\in \cX^\text{top}$. Let $v\in S$, such that $v < \xi$. Then, $\|g(x)-g(v)\| \leq L \| x-v\| \leq L \diam(\xi)$. This, along with \eqref{eq:beta}, shows that $g(x)\in B_{L\diam(\xi)}\left( E_{\Sigma(v)}(\mu(v),r(v))\right)$ and therefore, $(x,g(x))\in \tilde{G}_{L,r}$.
\end{proof}

We remark that Lemma~\ref{lem:approx} does not depend on the choice of cell complex $\cX$.
Thus we exploit the size of the geometric representation of cells to control the size of $G$.
Set 
\[
\varepsilon := \diam(\cX)\quad \text{and}\quad  \ell := \max\setof{\diam(E_{\Sigma(v)}(\mu(v),r(v)))\ |\ v\in S}.
\]

Given $\xi\in\cX^\text{top}$ let $\tilde{Q}(\xi) := \supp{\cF(\xi)} \cup Q(\xi) = \supp{\setof{\xi'\in\cX^\text{top}\mid |\xi'|\cap Q(\xi) \neq \emptyset }} \cup Q(\xi)$.

\begin{lem}
\label{lem:diamAxi}
Let $S \subset \cX^{(0)}$ be such that for each $\xi\in \cX^\text{top}$, there exists $v \in S$ such that $v < \xi$ and
let $\tilde{L}$ be chosen as in Lemma~\ref{lem:approx} and $L > \tilde{L}$. If $\xi\in\cX^\text{top}$, then
\[
\diam(\tilde{Q}(\xi)) < \ell + 2 L \varepsilon + 2 \varepsilon.
\]
\end{lem}

\begin{proof}
Let $\xi\in\cX^\text{top}$. Fix $v_0\in S$ such that $v_0<\xi$ and observe that, by equation \eqref{eq:Q},
\begin{align*}
\diam (Q(\xi))& = \diam \left( \bigcap_{v<\xi,\, v\in S} B_{L\diam{\xi}}\left(E_{\Sigma(v)}(\mu(v),r(v))\right)\right) \\
& \leq \diam\left( B_{L\diam{\xi}}\left(E_{\Sigma(v_0)}(\mu(v_0),r(v_0))\right)\right) \leq \ell + 2 L \varepsilon.
\end{align*}
The result follows from the fact that $\cF(\xi)$ is obtained by covering $Q(\xi) \cap X$ by elements of $\cX^\text{top}$.
\end{proof}

For the remaining of this section we assume the that $S = \cX^{(0)}$.

\begin{lem}
\label{lem:Gbound}
Let $\tilde{L}$ be chosen as in Lemma~\ref{lem:approx} and $L > \tilde{L}$, and let $\tilde{G} = \tilde{G}_{L, r}$ be defined as in equation \eqref{eq:defnGtilde}.
Then,
\[
\diam(\tilde{G}_x) < 2 (\ell + 2 L \varepsilon + 2 \varepsilon).
\]
\end{lem}

\begin{proof}
Note that $\tilde{G}_x = \bigcup \setof{ \{ x \} \times \tilde{Q}(\xi) \mid x \in \supp{\xi}, \xi \in \cX^\text{top}}$. If there is only one $\xi\in \cX^\text{top}$ such that $x \in \supp{\xi}$, then the result follows from Lemma~\ref{lem:diamAxi}. So assume that there are multiple $\xi \in \cX^\text{top}$ for which $x \in \supp{\xi}$. Let $(x, y), (x, y') \in \tilde{G}_x$. Then there exists $\xi,\xi' \in \cX^\text{top}$ such that $x \in \supp{\xi} \cap \supp{\xi'}$, $y \in \tilde{Q}(\xi)$, and $y' \in \tilde{Q}(\xi')$. Since $\supp{\xi} \cap \supp{\xi'} \neq \emptyset$, it follows that there exists $v_0 \in \cX^{(0)}$ such that $v_0 < \xi$ and $v_0 < \xi'$. Then, from the definition of $\tilde{Q}$ we have that $\mu(v_0) \in \tilde{Q}(\xi) \cap \tilde{Q}(\xi')$, and hence that $\tilde{Q}(\xi) \cap \tilde{Q}(\xi') \neq \emptyset$. Therefore it follows from Lemma~\ref{lem:diamAxi} that
\[
\|y -y'\| \leq \diam (\tilde{Q}(\xi)) + \diam (\tilde{Q}(\xi')) < 2 (\ell + 2 L \varepsilon + 2 \varepsilon),
\]
from which the result follows.
\end{proof}

Up to this point the construction is valid for any data set $\cT$ that satisfies assumption {\bf A}.
To control the size of the fibers of $\tilde{G}$ we recall \cite[Proposition 1]{freitas} that if a kernel $k$ of a $1$-dimensional Gaussian process is four times differentiable on the diagonal, and if we sample densely enough, then the posterior predictive variance $\sigma^2$ is uniformly bounded from above. 
More precisely, if the set of sample points is a $\gamma$-cover of $X$, then there exists a constant $Q^2\leq\sup_{X}{\partial ^2 _x\partial ^2 _y} k (x,y)|_{x=y}$ such that
\begin{equation}
\label{eq:upper_sigma}
\sup _{x\in X}\sigma(x)\leq \frac{Q\gamma^2}{4}.
\end{equation}

For $d$-dimensional outputs, by vectorizing the outputs and having a pre-specified kernel function $k$ for the Gaussian process, the prediction for an untried point can be obtained by a $d$-dimensional normal distribution with mean $\mu(x)$ and a covariance matrix $\Sigma(x)$. The maximum eigenvalue of $\Sigma(x)$ is bounded by ${\rm trace}(\Sigma(x))$ which equals to the summation of the variances in each dimension, so the result in \eqref{eq:upper_sigma} can be applied to each dimension. Therefore, given $\gamma >0$ and a set of sample points that is a $\gamma$-cover of $X$, there exists a constant $C= C(\gamma)$ such that
\begin{equation}
\label{eq:upper_sigma3}
\sup\left\{\lambda(x)\mbox{ a maximal eigenvalue of }\Sigma(x)\ |\ x\in X\right\}\leq  C (\gamma).
\end{equation}
Moreover,
\begin{equation}
\label{eq:gamma_conv}
\lim_{\gamma\to 0}C (\gamma)=0.
\end{equation}

\begin{proof}[Proof of Theorem~\ref{thm:convergence}]
Fix $L > \tilde{L}$, where $\tilde{L}$ satisfies the conditions in Lemma~\ref{lem:approx}.
Fix $\varepsilon_0 >0$ such that  
\begin{equation}
\label{eq:3le}
4(L+1)\varepsilon_0\leq\alpha/2.
\end{equation}
Consider $\cX$ a CW-structure on $X$ with $\diam(\cX)=\varepsilon<\varepsilon_0$ and a map $r:S\to(0,\infty)$ as in Lemma~\ref{lem:approx}.
Set $R :=\max\setof{r(v)\mid v\in S}$.
By \eqref{eq:gamma_conv} we can choose $\gamma_0 >0$ such that
\begin{equation}
\label{eq:cgamma}
C(\gamma_0)<\frac{\alpha ^2}{64 R}.
\end{equation}

By \cite[Theorem 3.7]{Alvarado21} there exists an $N_0\in\N$ such that any sample of size $N>N_0$ is a $\gamma_0$-cover of $X$ with probability greater than $1-\delta$. 
Thus, we assume that $N$, the number of data points in $\cT$, satisfies $N>N_0$.

Let $\tilde{G}$ be defined as in equation \eqref{eq:defnGtilde}. Note that the high probability inclusion \eqref{eq:incl} follows from Lemma~\ref{lem:approx}.
We conclude the proof by verifying $$\PP\left(\sup_{x\in X} \diam(\tilde{G}_x) < \alpha \right)>1-\delta.$$
With probability at least $1-\delta$, $\cT$ is a $\gamma_0$-cover of $X$.
Therefore, by Lemma~\ref{lem:Gbound} and \eqref{eq:upper_sigma3}
\begin{align*}
\diam(\tilde{G}_x) & \leq 2 \ell + 4 (L + 1) \varepsilon \leq 4 \sqrt{C(\gamma_0) R} + 4 (L +1) \varepsilon < \alpha
\end{align*}
where the last inequality follows from \eqref{eq:cgamma} and \eqref{eq:3le}.
\end{proof}

\section{Examples}
\label{sec:examples}
We conclude with examples demonstrating that our approach is capable of identifying  with high levels of confidence a broad range of dynamics from relatively few data points.

As discussed at the end of Section~\ref{sec:conley}, the Conley index computed via the multi-valued map $\cF$ is valid for any dynamical system for which $\cF$ is an outer approximation. Hence the results present in this section about fixed points, periodic orbits, connecting orbits, and chaotic dynamics occur for any dynamical system for which $\cF$ is an outer approximation, that is, the results of each example are valid for any sample path $h$ of the Gaussian process whose graph lies inside the region $G$ defined by equation \eqref{eq:defnG} in terms of $\cF$.

As indicated in Section~\ref{sec:conley}, the Conley index of $\cM\in\sM(\cF)$ in dimension $n$ is the rational canonical form of the linear map
\[
\cF_n \colon H_n \left(\down(\nu(\cM)),\down(\pred{\nu(\cM)});\F\right) \to H_n\left(\down(\nu(\cM)),\down(\pred{\nu(\cM)});\F\right).
\]
Since we are working with $1$-dimensional complexes $\cF_n = 0$ for all $n \geq 2$ and we can express the Conley index of $\cM$ as
\[
\Con_*(\cM; \F) \cong (p_0(x), p_1(x))
\]
where $p_n(x)$ is a monic polynomial \cite{bush:cowan:harker:mischaikow}. For the remainder of this discussion $\F = \Z_5$ (this allows us to distinguish whether the dynamics is orientation preserving or reversing).

We briefly mention a few standard results about the Conley index (see \cite{mischaikow:mrozek} for more details).
A \emph{trivial} Conley index in the $i$-th dimension takes the form $p_i(x) = x^k$ for some $k\in\N$.
\emph{However, to emphasize the triviality in the figures of this paper we write $p_i(x) =0$.}
If the Conley index is not trivial, i.e. $p_i(x) \neq x^k$, then the maximal invariant set in $\supp{\nu(\cM)}$ is nonempty (the converse is not true in general).
In the examples of this paper we make use of the following facts.
If the Conley index has the form $(x\pm 1, 0)$ or $(0, x\pm 1)$, then the maximal invariant set in $\supp{\nu(\cM)}$ contains a fixed point.
If the Conley index has the form $(x^T\pm 1,0)$ or $(0, x^T\pm 1)$, the the maximal invariant set in $\supp{\nu(\cM)}$ contains a periodic orbit of period $T$.
Let $\cM_i$ and $\cM_{i+1}$ be nodes in a Morse graph and assume that with respect to the poset order $\cM_{i+1}$ covers $\cM_{i}$. If
\begin{multline*}
\Con_*\left(\down(\nu(\cM_{i+1})),\down(\pred{\nu(\cM_i)});\F\right) \not\cong \\
\Con_*\left(\down(\nu(\cM_{i+1})),\down(\pred{\nu(\cM_{i+1})});\F\right)
\oplus
\Con_*\left(\down(\nu(\cM_{i})),\down(\pred{\nu(\cM_i)});\F\right)
\end{multline*}
then there exists a connecting orbit from $\cM_{i+1}$ to $\cM_{i}$.

As indicated above, the goal of this section is to show via examples that our approach can identify with high probability, interesting dynamics based on relatively few data points. To do this we pick a smooth function $f$ defined on an interval $X\subset \R$ that we know produces dynamics of interest, e.g., fixed points, periodic orbits, and chaotic dynamics. We then randomly (i.i.d.\ with respect to the uniform measure on $X$) generate data points $\cT = \setof{(x_n,f(x_n))\mid n=1,\ldots, N}$. In line with the fact that $f$ is smooth, we have chosen to use a smooth kernel, the squared exponential kernel $k(x,x') = \exp\left(-|x-x'|^2/\theta \right)$, and we \emph{assume} that the function $f$ is a sample path of the GP obtained from the data using MLE to estimate $\theta$.

For the examples in this section, we know the function $f$ from which the data is sampled. Therefore, we can verify that indeed $\sG(f) \subset G$. As a consequence, the dynamics that we report is valid for $f$. Of course, in applications we do not expect to be able to perform this explicit validation. We can only claim that if the Lipschitz constant $L$ used is large enough (see below) and if we \emph{assume} that $f$ is a sample path of the GP obtained from the data, then we have provided a lower bound on the confidence level that the dynamics identified via the homological calculations is valid for $f$. Without the assumption that $f$ is a sample path of the GP we can only claim that we have provided a lower bound on the probability that a dynamical system generated by a sample path of the obtained GP will exhibit the dynamics identified by our method.

Our method requires knowledge of the data set $\cT$ and the assumption that our choice of $L$ provides a bound for the Lipschitz constant of the unknown function $f$ with confidence $1 - \delta / 2$, that is, we assume that the probability in \eqref{eq:L} is at least $1 - \delta / 2$. We then select the set $S$ and construct the function $r$ such that the probability in equation \eqref{eq:conf} is also at least $1 - \delta / 2$. Then we construct the set $G$ in Theorem~\ref{thm:convergence} satisfying $\PP(\sG(g)\subset G) > 1 - \delta$ (see Lemma~\ref{lem:approx}). Hence the extracted dynamics has a confidence level of at least $1 - \delta$.

Although the function $f$ does not need to be specified for our method, in the examples we indicate $f$ to illustrate the fact that, at least in these cases, we can recover the dynamics of $f$ with very few data points. Since we do not know the values of $L_0$, $a$, and $b$ necessary for \eqref{eq:L}, in the examples we take $L$ to be at least twice the Lipschitz constant of $f$. For all the examples presented here we used $L=8$. The cell complex $\cX$ is obtained by considering a uniform decomposition of the interval $X$ into $2^B$ subintervals.

For the results presented here, our primary concern in the choice of $G$, which given a cell complex $\cX$ is equivalent to a choice of $\cF\colon \cX\mvmap \cX$, is to emphasize our ability to characterize dynamics with high levels of confidence. Thus, we impose a small $\delta$ in \eqref{eq:conf} on $S$. Unless otherwise stated, we use $\delta = 0.05$ and hence we obtain a $0.95$ (or $95\%$) confidence level.

\paragraph{Bistability}

To demonstrate that bistability can be identified via the Morse graph, we turn to Fig.~\ref{fig:2examples}~{\bf (a)} that was generated using $N=8$ data points sampled on the interval $X=[0,1]$ for the function $f(x) = 0.3 \arctan(8x - 4) + 0.5$.
The set $G$ was constructed using pointwise equal confidence intervals with $B = 9$ and $\delta=0.05$. Observe that the Morse graph has two minimal nodes $\cM(0)$ and $\cM(2)$. Under the poset isomorphism $\nu(\cM(0))$ and $\nu(\cM(2))$ are disjoint attracting blocks. Therefore the dynamics of any function $h$ with $\sG(h) \subset G$ exhibits at least bistability (it is possible that for a given $h$ there are additional attractors within $\nu(\cM(0))$ and/or $\nu(\cM(2))$.

The Conley indices for the bistability example are presented in Fig.~\ref{fig:2examples}~{\bf (a)}. The intervals defining the Morse sets are:
\[
|\nu(M_0)| = [0.09179688, 0.27929688],
\]
\[
|\nu(M_1)| = [0.76171875, 0.94140625],
\]
\[
|\nu(M_2)| = [0.75976562, 0.76171875],
\]
\[
|\nu(M_3)| = [0.50976562, 0.51171875],
\]
and
\begin{align*}
|\nu(M_4)| = & [0.49023438, 0.49218750] \cup [0.49414062, 0.49609375] \cup \\
& [0.49804688, 0.50390625] \cup [0.50585938, 0.50781250].
\end{align*}

\paragraph{Periodic Orbit}
As indicated in Section~\ref{sec:intro}, the Conley index can be used to identify periodic orbits of a given period. To demonstrate this, and to emphasize the importance of being able to choose $G$, we consider the logistic map $f(x) = 3.15 x (1-x)$. The global dynamics for $f$ is well understood. The points $x_0=0$ and $x_1 \approx 0.68$ are unstable fixed points, and all other initial conditions in $(0,1)$ limit to a stable period-2 orbit $\gamma$. 

\begin{figure}[!htbp]
\centering
\includegraphics[width=1.0\linewidth]{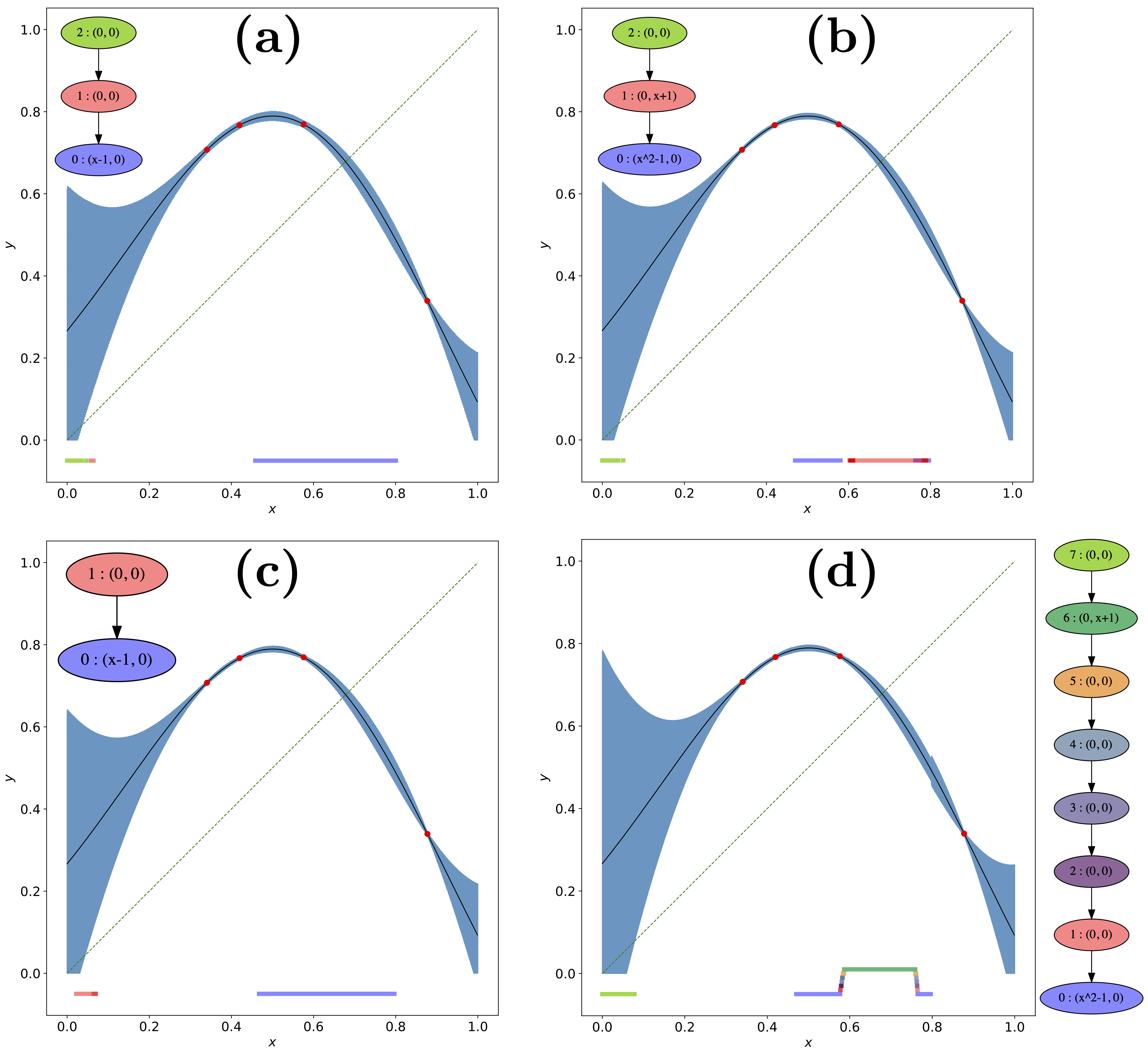}
\caption{In all figures elements of $\cT$ are indicated in red, the mean function $\mu$ is shown in black, and $G$ is shown in blue. The Morse graphs are indicated on the figures and the corresponding (color coded) regions  of phase space $\bar{\nu}(\cdot)$ are indicated at the bottom of each figure. {\bf (a)} The region $G$ is composed of squares of width $2^{-10} \length(X)$. The Conley index for $\cM_0$ identifies the existence of a fixed point with $95\%$ confidence. {\bf (b)} The region $G$ is composed of squares of width $2^{-11} \length(X)$. The Conley index for $\cM_0$ indicate the existence of a period two orbit with $95\%$ confidence. {\bf (c)} The region $G$ is composed of squares of width $2^{-11} \length(X)$ and in this case the Conley index for $\cM_0$ identifies the existence of a fixed point with $97.5\%$ confidence, but we do not detect the periodic orbit. {\bf (d)} The region $G$ is composed of squares of width $2^{-11} \length(X)$. The Conley index for $\cM_0$ indicate the existence of a period two orbit with $97.5\%$ confidence.}
\label{fig:ex2}
\end{figure}

To apply our techniques, set $X=[0, 1]$, $N=4$, $B=10$, and $\delta =0.05$ (see Fig.~\ref{fig:ex2}~{\bf (a)}). Choosing pointwise equal confidence intervals for this value of $\delta$ leads to a Morse graph $\sM(\cF) = \{\cM_0,\cM_1, \cM_2 \mid \cM_0< \cM_1 < \cM_2 \}$ where 
$x_0 \in \supp{\nu(M_2)}$ and $\supp{\nu(M_0)}$ contains both $x_1$ and $\gamma$.
The Conley index of $\cM_0$ identifies the existence of a fixed point. 
Since $\delta =0.05$, this description of the dynamics is valid
with at least $95\%$ confidence.

What is missing from this description is the identification of a period-2 orbit.
Motivated by Theorem~\ref{thm:convergence}, we repeat the computations with $B = 11$ and detected the existence of a periodic orbit with at least $95\%$ confidence (see Fig.~\ref{fig:ex2}~{\bf (b)}).
We can also become more ambitious and seek a confidence level of $97.5\%$, i.e., setting $\delta = 0.025$.
For $B = 11$ we fail to identify the periodic orbit (see Fig.~\ref{fig:ex2}~{\bf (c)}).
Rather than increasing the subdivision, we choose nonuniform confidence intervals; where we make the images of $\cF(\xi)$ smaller if $\supp{\xi}\subset \supp{\nu(M_0)}$ and larger elsewhere (see Fig.~\ref{fig:ex2}~{\bf (d)}). 
The resulting Morse graph contains eight nodes where $x_0\in \supp{\nu(\cM_7)}$, $x_1 \in \supp{\nu(\cM_6)}$, and the Conley index of $\cM_0$ indicates the existence of a period-2 orbit. 

The Conley indices for the period-$2$ orbit example are presented in Fig.~\ref{fig:ex2}. The intervals defining the Morse sets are: For Fig.~\ref{fig:ex2}~{\bf (a)}:
\[
|\nu(M_0)| = [0.45800781, 0.80078125],
\]
\[
|\nu(M_1)| = [0.03613281, 0.03710938],
\]
and
\[
|\nu(M_2)| = [0.00000000, 0.03320312] \cup [0.03417969, 0.03515625];
\]
for Fig.~\ref{fig:ex2}~{\bf (b)}:
\[
|\nu(M_0)| = [0.46972656, 0.58056641] \cup [0.76269531, 0.79638672],
\]
\begin{align*}
|\nu(M_1)| = & [0.58056641, 0.58105469] \cup [0.58154297, 0.58203125] \cup [0.58251953, 0.58300781] \cup \\
& [0.58349609, 0.58398438] \cup [0.58447266, 0.58496094] \cup [0.58544922, 0.76074219] \cup \\
& [0.76123047, 0.76171875] \cup [0.76220703, 0.76269531],
\end{align*}
and
\[
|\nu(M_2)| = [0.00000000, 0.03808594] \cup [0.03857422, 0.03906250];
\]
for Fig.~\ref{fig:ex2}~{\bf (c)}:
\[
|\nu(M_0)| = [0.46728516, 0.79687500],
\]
and
\[
|\nu(M_1)| = [0.00000000, 0.04199219] \cup [0.04248047, 0.04296875];
\]
for Fig.~\ref{fig:ex2}~{\bf (d)}:
\[
|\nu(M_0)| = [0.47021484, 0.57666016] \cup [0.76464844, 0.79638672],
\]
\[
|\nu(M_1)| = [0.57666016, 0.57714844] \cup [0.57763672, 0.57812500] \cup [0.76416016, 0.76464844],
\]
\[
|\nu(M_2)| = [0.57861328, 0.57910156] \cup [0.57958984, 0.58007812] \cup [0.76318359, 0.76367188],
\]
\[
|\nu(M_3)| = [0.58056641, 0.58105469] \cup [0.76220703, 0.76269531],
\]
\[
|\nu(M_4)| = [0.58154297, 0.58203125] \cup [0.58251953, 0.58300781] \cup [0.76123047, 0.76171875],
\]
\begin{align*}
|\nu(M_5)| = & [0.58349609, 0.58398438] \cup [0.58447266, 0.58496094] \cup [0.58544922, 0.58593750] \cup \\
& [0.75927734, 0.75976562] \cup [0.76025391, 0.76074219],
\end{align*}
\begin{align*}
|\nu(M_6)| = & [0.58642578, 0.58691406] \cup [0.58740234, 0.58789062] \cup [0.58837891, 0.58886719] \cup \\
& [0.58935547, 0.58984375] \cup [0.59033203, 0.59082031] \cup [0.59130859, 0.75683594] \cup \\
& [0.75732422, 0.75781250] \cup [0.75830078, 0.75878906],
\end{align*}
and
\[
|\nu(M_7)| = [0.00000000, 0.07812500].
\]

\begin{figure}[!htbp]
\centering
\includegraphics[width=0.6\linewidth]{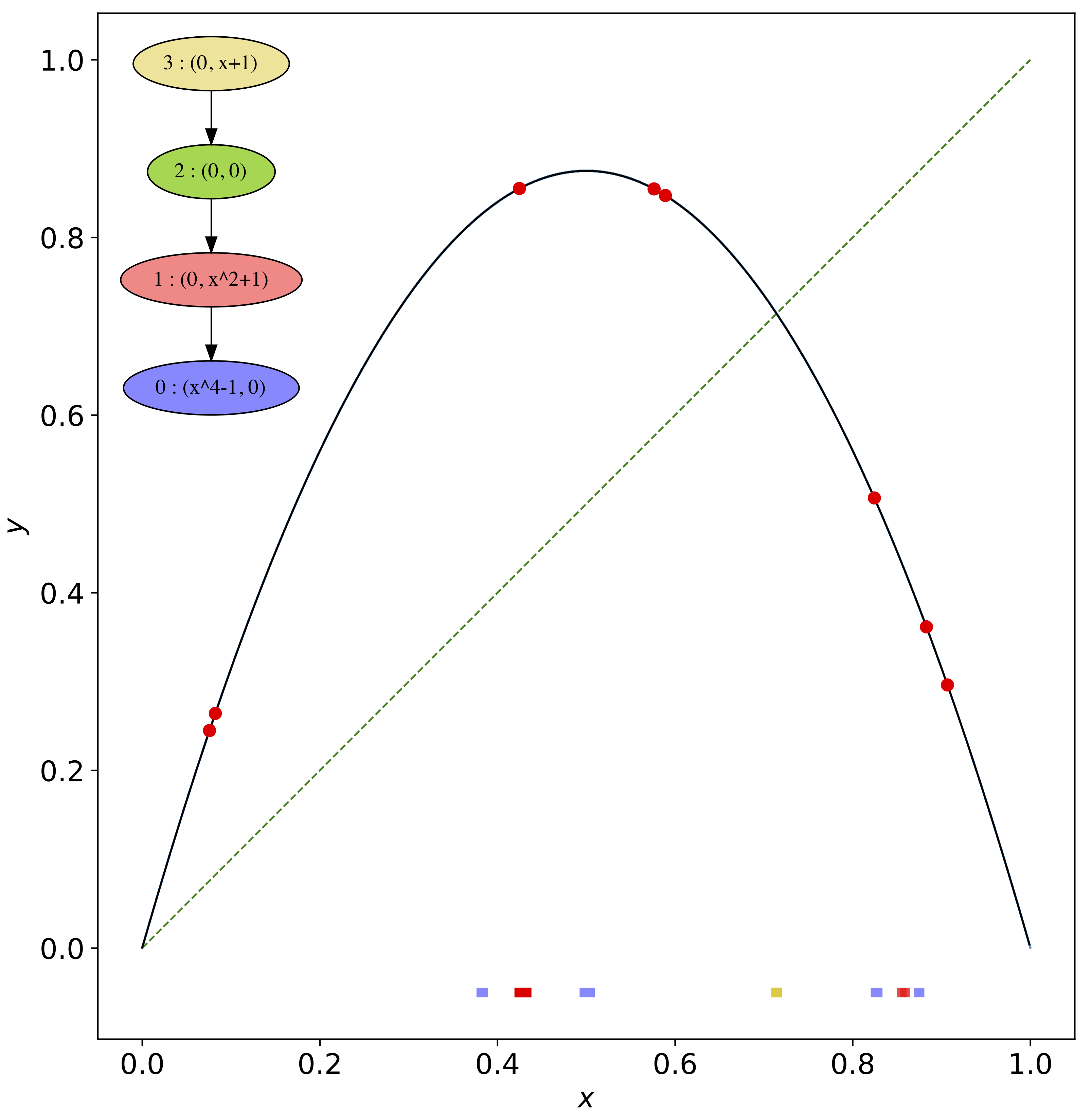}
\caption{In both figures elements of $\cT$ are indicated in red, the mean function $\mu$ is shown in black and $G$ is shown in blue where the width of $G$ is approximately similar to the graphical representation of $\mu$. The Morse graphs are indicated on the figures and the corresponding (color coded) regions  of phase space $\bar{\nu}(\cdot)$ are indicated at the bottom of each figure. The region $G$ is composed of squares of width $2^{-15} \length(X)$. The Conley indices for $\cM_1$ and $\cM_0$ indicated the existence of period 2 and period 4 orbits, respectively, with $95\%$ confidence.}
\label{fig:ex3}
\end{figure}

\paragraph{Connecting Orbits}
Consider again the logistic map $f(x) = 3.5 x (1-x)$, $X=[0, 1]$, $N = 8$, $B = 15$, and $\delta =0.05$. Using pointwise equal confidence intervals, we obtain the Morse graph in Fig.~\ref{fig:ex3}. The Conley indices for $M_1$ and $M_0$ indicated the existence of period-2 and period-4 orbits, respectively, with at least $95\%$ confidence. The Conley index also indicates the existence of a fixed point for $M_3$ and a connecting orbit from $M_3$ to $M_1$.

The Conley indices for the period-$4$ orbit and connecting orbit example are presented in Fig.~\ref{fig:ex3}. The intervals defining the Morse sets are:
\begin{align*}
|\nu(M_0)| = & [0.38165283, 0.38403320] \cup [0.49774170, 0.50427246] \cup \\
& [0.82553101, 0.82818604] \cup [0.87463379, 0.87533569],
\end{align*}
\begin{align*}
|\nu(M_1)| = & [0.42437744, 0.42440796] \cup [0.42443848, 0.42446899] \cup [0.42449951, 0.42453003] \cup \\
& [0.42456055, 0.42459106] \cup [0.42462158, 0.42465210] \cup [0.42468262, 0.42471313] \cup \\
& [0.42474365, 0.42477417] \cup [0.42480469, 0.43246460] \cup [0.43249512, 0.43252563] \cup \\
& [0.43255615, 0.43258667] \cup [0.43261719, 0.43264771] \cup [0.43267822, 0.43270874] \cup \\
& [0.43273926, 0.43276978] \cup [0.43280029, 0.43283081] \cup [0.43286133, 0.43289185] \cup \\
& [0.43292236, 0.43295288] \cup [0.85528564, 0.85531616] \cup [0.85534668, 0.85885620] \cup \\
& [0.85888672, 0.85891724],
\end{align*}
\[
|\nu(M_2)| = [0.71362305, 0.71365356] \cup [0.71490479, 0.71493530],
\]
and
\begin{align*}
|\nu(M_3)| = & [0.71368408, 0.71371460] \cup [0.71374512, 0.71377563] \cup [0.71380615, 0.71475220] \cup \\
& [0.71478271, 0.71481323] \cup [0.71484375, 0.71487427].
\end{align*}

The direct sum of the Conley indices of $\nu(\cM_3)$ and  $\nu(\cM_1)$ 
\[
\Con_*\left(\down(\nu(\cM_{3})),\down(\pred{\nu(\cM_{3})});\F\right)
\oplus
\Con_*\left(\down(\nu(\cM_{1})),\down(\pred{\nu(\cM_1)});\F\right)
\]
can be represented by
\[
{\small
\setcounter{MaxMatrixCols}{3}
\begin{bmatrix}
-1&     0&     0\\
0&     0&    -1\\
0&     1&     0\\
\end{bmatrix}.
}
\]
The Conley index of the connecting orbit set
\[
\Con_*\left(\down(\nu(\cM_{3})),\down(\pred{\nu(\cM_1)});\F\right) 
\]
can be represented by
\[
{\small
\setcounter{MaxMatrixCols}{14}
\begin{bmatrix}
-1&     0&     0&     0&     0&     0&     0&     0&     0&     0&  0&  0&  0&  0\\
 1&     0&     0&     0&     0&     0&     0&     0&     0&     0&  0&  0&  0&  0\\
 0&    -1&     0&     0&     0&     0&     0&     0&     0&     0&  0&  0&  0&  0\\
 0&     0&    -1&     0&     0&     0&     0&     0&     0&     0&  0&  0&  0&  0\\
 0&     0&     0&    -1&     0&     0&     0&     0&     0&     0&  0&  0&  0&  0\\
 0&     0&     0&     0&    -1&     0&     0&     0&     0&     0&  0&  0&  0&  0\\
 0&     0&     0&     0&     0&     1&     0&     0&     0&     0&  0&  0&  0&  0\\
 0&     0&     0&     0&     0&     0&     1&     0&     0&     0&  0&  0&  0&  0\\
 0&     0&     0&     0&     0&     0&     0&     1&     0&     0&  0&  0&  0&  0\\
 0&     0&     0&     0&     0&     0&     0&     0&     1&     0&  0&  0&  0&  0\\
 0&     0&     0&     0&     0&     0&     0&     0&     0&    -1&  0&  0&  0&  0\\
 0&     0&     0&     0&     0&     0&     0&     0&     0&     0&  1&  0&  0&  0\\
 0&     0&     0&     0&     0&     0&     0&     0&     0&     0&  0&  1&  0&  1\\
 0&     0&     0&     0&     0&     0&     0&     0&     0&     0&  0&  0& -1&  0\\
\end{bmatrix}.
}
\]
From these computations
\begin{multline*}
\Con_*\left(\down(\nu(\cM_{3})),\down(\pred{\nu(\cM_1)});\F\right) \not\cong \\
\Con_*\left(\down(\nu(\cM_{3})),\down(\pred{\nu(\cM_{3})});\F\right)
\oplus
\Con_*\left(\down(\nu(\cM_{1})),\down(\pred{\nu(\cM_1)});\F\right)
\end{multline*}
which indicates the existence of a connecting orbit from the Morse set in $|\nu(\cM_3)|$ to the Morse set in $|\nu(\cM_1)|$.

\paragraph{Chaotic dynamics}
Consider the function $f(x) = 2 e^{-5 (x - 1)^2}$, $X=[-0.2, 2.3]$, $N=10$, $B=10$, and $\delta = 0.05$. Using pointwise equal confidence intervals, as indicated in Fig.~\ref{fig:2examples}~{\bf (b)}, we obtain a Morse graph with five nodes.
While the Conley index of $\cM(4)$ is trivial, it consists of multiple disjoint intervals and the index map $\cF_*\colon H_*(\nu(\cM(4)),\pred{\nu(\cM(4))};\Z_5)\to H_*(\nu(\cM(4)),\pred{\nu(\cM(4))};\Z_5)$ can be used to capture the chaotic dynamics (c.f.\ \cite{day} and references therein).

The Conley indices for the chaotic dynamics example are presented in Fig.~\ref{fig:2examples}~{\bf (b)}. In this case the Conley index of $\nu(\cM(4))$ is trivial, however the index map 
\[
\cF_1\colon H_1\left(\nu(\cM(4)),\pred{\nu(\cM(4))};\Z_5\right)\to H_1(\nu(\cM(4)),\pred{\nu(\cM(4))};\Z_5)\] 
represented by
\[
{\small
\setcounter{MaxMatrixCols}{20}
\begin{bmatrix}
0 &  1 &  0 &  0 &  0 &  0 &  0 &  0 &  0 &  0 & -1 &  0 &  0 &  0 &  0 &  0 &  0 &  0 &  0 &  0 \\
0 &  1 &  0 &  0 &  0 &  0 &  0 &  0 &  0 &  0 & -1 &  0 &  0 &  0 &  0 &  0 &  0 &  0 &  0 &  0 \\
0 &  0 &  0 &  0 &  0 &  0 &  0 &  0 &  0 &  0 &  0 &  0 &  0 &  0 &  0 &  0 &  0 &  1 &  0 & -1 \\
0 &  1 &  0 &  0 &  0 &  0 &  0 &  0 &  0 &  0 & -1 &  0 &  0 &  0 &  0 &  0 &  0 &  0 &  0 &  0 \\
0 &  0 &  0 &  0 &  0 &  0 &  0 &  0 &  0 &  0 &  0 &  0 &  0 &  0 &  0 &  0 &  0 &  1 &  0 & -1 \\
0 &  0 &  0 &  1 &  0 &  0 &  0 &  0 &  0 & -1 &  0 &  0 &  0 &  0 &  0 &  0 &  0 &  0 &  0 &  0 \\
0 &  0 &  0 &  0 &  0 &  0 &  0 &  0 &  0 &  0 &  0 &  0 &  0 &  0 &  0 &  0 &  0 &  1 &  0 & -1 \\
0 &  0 &  0 &  1 &  0 &  0 &  0 &  0 &  0 & -1 &  0 &  0 &  0 &  0 &  0 &  0 &  0 &  0 &  0 &  0 \\
0 &  0 &  0 &  0 &  0 &  0 &  0 &  0 &  0 &  0 &  0 &  0 &  0 &  0 &  0 &  0 &  0 &  1 &  0 & -1 \\
0 &  0 &  0 &  0 &  0 &  0 &  0 &  0 &  0 &  0 &  0 &  0 &  0 &  0 &  0 &  0 &  0 &  1 &  0 & -1 \\
0 &  0 &  0 &  0 &  0 &  0 &  0 &  0 &  0 &  0 &  0 &  0 &  0 &  0 &  0 &  0 &  0 &  1 &  0 & -1 \\
0 &  0 &  0 &  0 &  0 &  0 &  0 &  0 &  0 &  0 &  0 &  0 &  0 &  0 &  0 &  0 &  0 &  1 &  0 & -1 \\
0 &  0 &  0 &  0 & -1 &  0 &  0 &  1 &  0 &  0 &  0 &  0 &  0 &  0 &  0 &  0 &  0 &  0 &  0 &  0 \\
0 &  0 &  0 &  0 & -1 &  0 &  0 &  1 &  0 &  0 &  0 &  0 &  0 &  0 &  0 &  0 &  0 &  0 &  0 &  0 \\
0 &  0 &  0 &  0 & -1 &  0 &  0 &  1 &  0 &  0 &  0 &  0 &  0 &  0 &  0 &  0 &  0 &  0 &  0 &  0 \\
0 &  0 &  0 &  0 &  0 &  0 &  0 &  0 &  0 &  0 &  0 &  0 &  0 &  1 &  0 &  0 &  0 &  0 &  0 & -1 \\
0 &  0 &  0 &  0 & -1 &  0 &  0 &  1 &  0 &  0 &  0 &  0 &  0 &  0 &  0 &  0 &  0 &  0 &  0 &  0 \\
0 &  0 &  0 &  0 & -1 &  0 &  0 &  1 &  0 &  0 &  0 &  0 &  0 &  0 &  0 &  0 &  0 &  0 &  0 &  0 \\
0 &  0 &  0 &  0 & -1 &  0 &  0 &  1 &  0 &  0 &  0 &  0 &  0 &  0 &  0 &  0 &  0 &  0 &  0 &  0 \\
0 &  0 &  0 &  0 &  0 &  0 &  0 &  0 &  0 &  0 &  0 &  0 &  0 &  1 &  0 &  0 &  0 &  0 &  0 & -1 \\
\end{bmatrix}
}
\]
indicates the existence of chaotic dynamics (see \cite{day} and references therein).

The intervals defining the Morse sets are:
\[
|\nu(M_0)| = [-0.01933594, 0.10517578],
\]
\[
|\nu(M_1)| = [0.10517578, 0.10761719],
\]
\[
|\nu(M_2)| = [0.11005859, 0.11250000],
\]
\[
|\nu(M_3)| = [0.43232422, 0.43476562],
\]
and
\begin{align*}
|\nu(M_4)| = & [0.43720703, 0.43964844] \cup [0.44208984, 0.47871094] \cup [0.48115234, 0.50312500] \cup \\
& [0.52509766, 0.52753906] \cup [0.52998047, 0.57148437] \cup [0.67158203, 0.67402344] \cup \\
& [0.67646484, 0.71796875] \cup [0.72041016, 0.72285156] \cup [0.73017578, 0.73261719] \cup \\
& [0.73505859, 0.78144531] \cup [0.78388672, 0.78632812] \cup [1.19404297, 1.19648437] \cup \\
& [1.19892578, 1.33320312] \cup [1.42353516, 1.42597656] \cup [1.42841797, 1.46992187] \cup \\
& [1.47236328, 1.47480469] \cup [1.49189453, 1.49433594] \cup [1.49677734, 1.51875000] \cup \\
& [1.52119141, 1.55781250] \cup [1.56025391, 1.56269531].
\end{align*}

\section{Construction of $G$ for examples in Section~\ref{sec:examples}}

\begin{figure}[!htbp]
\begin{center}
\begin{tikzpicture}
[scale=0.3]


\fill[blue!50] (1,9) rectangle (5,29);
\fill[blue!50] (9,13) rectangle (13,33);
\fill[blue!50] (17,13) rectangle (21,37);
\fill[blue!50] (25,17) rectangle (29,33);

\fill[teal!50] (5,1) rectangle (9,41);
\fill[teal!50] (13,5) rectangle (17,45);
\fill[teal!50] (21,5) rectangle (25,41);

\draw[dotted, thick] (1,0) -- (1,46);
\draw(1,-0.7) node{\small{$v_{2i-3}$}};
\draw[dotted, thick] (5,0) -- (5,46);
\draw(5,-0.7) node{\small{$v_{2i-2}$}};
\draw[dotted, thick] (9,0) -- (9,46);
\draw(9,-0.7) node{\small{$v_{2i-1}$}};
\draw[dotted, thick] (13,0) -- (13,46);
\draw(13,-0.7) node{\small{$v_{2i}$}};
\draw[dotted, thick] (17,0) -- (17,46);
\draw(17,-0.7) node{\small{$v_{2i+1}$}};
\draw[dotted, thick] (21,0) -- (21,46);
\draw(21,-0.7) node{\small{$v_{2i+2}$}};
\draw[dotted, thick] (25,0) -- (25,46);
\draw(25,-0.7) node{\small{$v_{2i+3}$}};
\draw[dotted, thick] (29,0) -- (29,46);
\draw(29,-0.7) node{\small{$v_{2i+4}$}};

\draw[dotted, thick] (0,1) -- (30,1);
\draw[dotted, thick] (0,5) -- (30,5);
\draw[dotted, thick] (0,9) -- (30,9);
\draw[dotted, thick] (0,13) -- (30,13);
\draw[dotted, thick] (0,17) -- (30,17);
\draw[dotted, thick] (0,21) -- (30,21);
\draw[dotted, thick] (0,25) -- (30,25);
\draw[dotted, thick] (0,29) -- (30,29);
\draw[dotted, thick] (0,33) -- (30,33);
\draw[dotted, thick] (0,37) -- (30,37);
\draw[dotted, thick] (0,41) -- (30,41);
\draw[dotted, thick] (0,45) -- (30,45);

\draw[dashed, thick] (0,16) .. controls (10,26) and (20,26) .. (30,24);
\draw (31,24) node{$\mu(x)$};

\draw[thick] (27,17.5) -- (27,31.5);
\draw[thick] (19,17) -- (19,33);
\draw[thick] (11,14) -- (11,32);
\draw[thick] (3,10) -- (3,28);

\draw[red, thick] (3,10) -- (8,0);
\draw[red, thick] (11,14) -- (4,0);
\draw[red, thick] (11,14) -- (17,2);
\draw[red, thick] (19,17) -- (12,3);
\draw[red, thick] (19,17) -- (27,1);
\draw[red, thick] (27,17.5) -- (22,7.5);

\draw[red, thick] (3,28) -- (9,40);
\draw[red, thick] (11,32) -- (6,42);
\draw[red, thick] (11,32) -- (15,44);
\draw[red, thick] (19,33) -- (13,45);
\draw[red, thick] (19,33) -- (25,45);
\draw[red, thick] (27,31.5) -- (22,41.5);
\end{tikzpicture}
\end{center}
\caption{Construction of $\cF$. 
The dashed line represents $\mu(x)$ derived from the data $\cT$. 
The dotted lines indicate the cell complex $\cX \times \cX$.
The solid dark lines indicate the confidence intervals at the midpoints for the odd numbered intervals.
The blue regions indicate the value of $\cF$ on odd numbered edges. 
The red rays have slope $\pm L$ and represents high probability Lipschitz bounds on sample paths.
The teal regions indicate the value of $\cF$ on even numbered edges.}
\label{fig:MVMAP}
\end{figure}
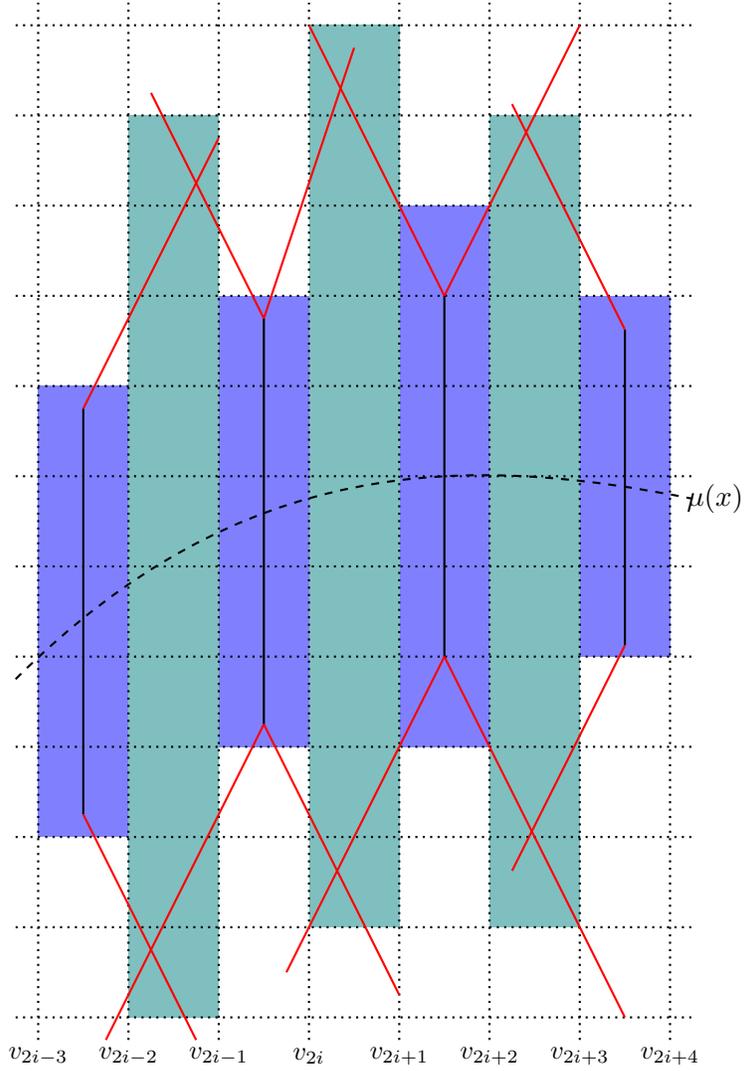

For the sake of clarity, the discussion of Conley theory in Section~\ref{sec:conley} and the probabilistic bounds in Section~\ref{sec:surrogate} avoided options that can improve computational efficacy.
Turning to the details of the computations for the results reported in Section~\ref{sec:examples}, we take explicit advantage of some of these options, and we exploit the fact that we are working with one-dimensional dynamics to pictorially explain the computational enhancements. 
As a starting point, see Fig.~\ref{fig:MVMAP} where the dashed curve indicates $\mu$.

Throughout this discussion $X=[\alpha,\beta]\subset \R$.
To avoid cumbersome notation we identify the cellular complex $\cX$ with a uniform discretization of $[\alpha,\beta]$ into $2^B$ subintervals.
We denote the vertices  of $\cX$ by $\cX^{(0)} = \setof{v_i\mid i=0,\ldots, 2^B}$ and the edges  by $$\cX^{(1)} = \setof{e_i = [v_i,v_{i+1}]\mid i=0,\ldots, 2^B-1}.$$
The set of midpoints of the odd intervals is denoted by $$S=\setof{m_{2i+1} = (v_{2i+1} + v_{2i+2})/2\mid i=0,\ldots, 2^{B-1}-1}.$$
The length of each interval is $\epsilon := 2^{-B}(\beta-\alpha)$.
The cubical grid on a subportion of $[\alpha,\beta]\times [\alpha,\beta]$ is indicated via the dotted lines in Fig.~\ref{fig:MVMAP}.

Fix $\delta \in (0,1)$.
As discussed in Section~\ref{sec:GP}, choose $r\colon S \to (0,\infty)$ such that equation \eqref{eq:conf} is satisfied for all $m_i \in S$.
For notational convenience we set
\[
\underline{w}_{2i+1} := \mu(m_{2i+1})-r(m_{2i+1})\sigma(m_{2i+1})
\quad\text{and}\quad
\overline{w}_{2i+1} := \mu(m_{2i+1})+r(m_{2i+1})\sigma(m_{2i+1}).
\]

We construct a multivalued map $\cF\colon \cX^\text{top} \mvmap \cX^\text{top}$ in two steps. Note that $\cX^\text{top} = \cX^{(1)}$ in the settings of this section.

\smallskip
\noindent
\emph{Step 1.} Let $\cX_\text{odd}^{(1)} = \setof{e_{2i+1}\mid i=0,\ldots, 2^{B-1}-1}\subset \cX^{(1)}$.
For each $e_{2i+1}\in \cX_\text{odd}^{(1)}$ define 
\[
\cF(e_{2i+1}) = \setof{e\in \cX^{(1)} \mid e \cap [\underline{w}_{2i+1},\overline{w}_{2i+1}] \neq \emptyset
}.
\]
In Fig.~\ref{fig:MVMAP} the black lines are used to indicate $[\underline{w}_{2i+1},\overline{w}_{2i+1}]$ and the light blue shaded regions designate $\cF(e_{2i+1})$.

\smallskip
\noindent
\emph{Step 2.} Observe that we are guaranteed with probability at least $1-\delta$, that $g(m_{2i+1})$ will lie in the black lines and hence blue regions of Fig.~\ref{fig:MVMAP}.
To gain control of $\sG(g)$ over $X \setminus S$ make use of the bound given by equation \eqref{eq:L}.
As stated in Section~\ref{sec:examples}, we assume that $L$ is large enough so that the probability given in equation \eqref{eq:L} is at least $(1-\delta)^{1/2}$. We also assume that $L$ is large enough so that the rays defined below intersect.
For each point $m_{2i+1}$ we consider the four rays
\begin{equation}
\label{eq:rays1}
\overline{h}_{2i+1}^+(s) = \left(m_{2i+1},\overline{w}_{2i+1}\right) + s(1,L)\quad\text{and}\quad 
\overline{h}_{2i+1}^-(s) = \left(m_{2i+1},\overline{w}_{2i+1}\right) + s(-1,L)
\end{equation}
and
\begin{equation}
\label{eq:rays2}
\underline{h}_{2i+1}^+(s) = \left(m_{2i+1},\underline{w}_{2i+1}\right) + s(1,-L)\quad\text{and}\quad 
\underline{h}_{2i+1}^-(s) = \left(m_{2i+1},\underline{w}_{2i+1}\right) + s(-1,-L)
\end{equation}
parameterized by $s\geq 0$.
These are shown as red lines in Fig.~\ref{fig:MVMAP}.

Observe that the rays $\overline{h}_{2i-1}^+$ and $\overline{h}_{2i+1}^-$ intersect at
\begin{align*}
\overline{h}_{2i-1}^+\left(\epsilon + \frac{\overline{w}_{2i+1}-\overline{w}_{2i-1}}{2L}\right) =
\overline{h}_{2i+1}^-\left(\epsilon - \frac{\overline{w}_{2i+1}-\overline{w}_{2i-1}}{2L}\right) = \\
\left(
m_{2i-1}+ \epsilon + \frac{\overline{w}_{2i+1}-\overline{w}_{2i-1}}{2L},
\overline{w}_{2i-1} + \epsilon L + \frac{\overline{w}_{2i+1}-\overline{w}_{2i-1}}{2}
\right)
\end{align*}
and the rays $\underline{h}_{2i-1}^+$ and $\underline{h}_{2i+1}^-$ intersect at
\begin{align*}
\underline{h}_{2i-1}^+\left(\epsilon -\frac{\underline{w}_{2i+1}-\underline{w}_{2i-1}}{2L}\right) = 
\underline{h}_{2i+1}^-\left(\epsilon +\frac{\underline{w}_{2i+1}-\underline{w}_{2i-1}}{2L}\right) = \\
\left(
m_{2i-1}+ \epsilon -\frac{\underline{w}_{2i+1}-\underline{w}_{2i-1}}{2L},
\underline{w}_{2i-1} - \epsilon L + \frac{\underline{w}_{2i+1}-\underline{w}_{2i-1}}{2}
\right)
\end{align*}

For $i > 0$ define
\[
Q(e_{2i}) = \left[\underline{w}_{2i-1} - \epsilon L + \frac{\underline{w}_{2i+1}-\underline{w}_{2i-1}}{2}, \overline{w}_{2i-1} +\epsilon L + \frac{\overline{w}_{2i+1}-\overline{w}_{2i-1}}{2}\right],
\]
and
\[
\cF(e_{2i}) = \setof{e\in \cX^{(1)} \mid e \cap Q(e_{2i}) \neq \emptyset}.
\]

For $i = 0$ define
\[
Q(e_{0}) = \left[\underline{w}_{1} - \frac{3}{2} \epsilon L, \overline{w}_{1} + \frac{3}{2} \epsilon L \right],
\]
and
\[
\cF(e_{0}) = \setof{e\in \cX^{(1)} \mid e \cap
Q(e_{0}) \neq \emptyset},
\]
where $\underline{w}_{1} - 3/2 \epsilon L$ and $\overline{w}_{1} + 3/2 \epsilon L$ are the second components of the intersections of the line $x = v_0$ with the rays $\underline{h}_{1}^-(s)$ and $\overline{h}_{1}^-(s)$, respectively.
The teal regions in Fig.~\ref{fig:MVMAP} indicate $\cF(e_{2i})$.

\smallskip
\noindent\emph{Remark.}
As a consequence of Steps 1 and 2 we have defined the acyclic multivalued map $\cF\colon \cX^\text{top} \mvmap \cX^\text{top}$ that is used to identify the Morse graphs, lattices of attractors, and compute Conley indices.

\smallskip
In the spirit of Section~\ref{sec:conley} of the main text, define
\[
G = \bigcup_{e \in \cX^{(1)}} e \times \cF(e) \subset X \times X,
\]
and
\[
\tilde{G} = G \cup \left( \bigcup_{e \in \cX^{(1)}} e \times Q(e) \right) \subset X \times \R,
\]
where $Q(e)$ is defined in \emph{Step 2} for $e \in \cX_\text{even}^{(1)}$ and for $e_{2i+1} \in \cX_\text{odd}^{(1)}$ we define
\[
Q(e_{2i+1}) := \left[ \underline{w}_{2i+1} - \frac{1}{2} \epsilon L, \overline{w}_{2i+1} + \frac{1}{2} \epsilon L \right],
\]
where $\underline{w}_{2i+1} - \frac{1}{2} \epsilon L$ and $\overline{w}_{2i+1} + \frac{1}{2} \epsilon L$ are the second components of the intersections of the rays with the lines $x = v_{2 i + 1}$ and $x = v_{2 i + 1}$, respectively.

Observe that by equations \eqref{eq:conf} and \eqref{eq:L} and an argument analogous to the proof of Lemma~\ref{lem:approx}
\[
\PP(\sG(g) \subset \tilde{G})> 1-\delta.
\]

As can be seen from Fig.~\ref{fig:MVMAP}, $\cF$ is not an outer approximation for every $g$ such that $\sG(g) \subset \tilde{G}$ (this could be achieved by enlarging the images of $\cF(e_{2i+1})$, but at the risk of losing information about the structure of the dynamics).
However, we note that if $(x,y) \in \tilde{G}$, then $(x,y)$ is within distance $\epsilon / 2$ from $G$.
Therefore, we can apply the results of \cite[Section 5]{batko:mischaikow:mrozek:przybylski} to conclude that the Conley index implications about the dynamics computed using $\cF$ are valid for a sample path $h \colon X \to X$ if $\sG(h) \subset \tilde{G}$.

\section{Concluding Remarks}
\label{sec:conclusions}
As is indicated in Section~\ref{sec:examples}, our framework is capable of identifying the fundamental building blocks of traditional nonlinear dynamical systems with high levels of confidence based on few data points.
The most obvious criticism is that we restricted our examples to one-dimensional dynamics.
This was for the sake of clarity; the results described in Sections~\ref{sec:GP}-\ref{sec:surrogate} are dimension independent.
Computations of the type described in Section~\ref{sec:conley}, using cubical complexes can be done routinely in systems of dimension four or less \cite{vieira}.
They have also been used for the rigorous analysis of infinite-dimensional systems \cite{day:junge:mischaikow} indicating that, at least conceptually, it is the intrinsic, as opposed to extrinsic dimension, of the dynamics that determines computability.

There are a variety of closely related open problems that arise from our approach. 
The geometry of isolating blocks in one-dimension is reasonably simple, a finite collection of closed intervals.
In higher dimensions the geometry can be much more complicated, which raises the question of estimates relating the dynamics, the number of data points, and uncertainty bounds.
Even heuristics for optimal sampling methods to identify attractor block lattices is not obvious.

\section*{Data Availability}
The code to perform the computations and generate the figures is available at \url{https://github.com/marciogameiro/GP_MorseGraph}.

\section*{Acknowledgments}
The work of B.B. was partially supported by DARPA contract HR0011-16-2-0033 and by the Polish National Science Center under Opus Grant No. 2019/35/B/ST1/00874. The work of M.G., Y.H., E.V., and K.M. was partially supported by the National Science Foundation under awards DMS-1839294 and HDR TRIPODS award CCF-1934924, DARPA contract HR0011-16-2-0033, National Institutes of Health award R01 GM126555, and Air Force Office of Scientific Research under award numbers FA9550-23-1-0011 and AWD00010853-MOD002. M.G. was also supported by FAPESP grant 2019/06249-7 and CNPq grant 309073/2019-7. K.M. was also supported by a grant from the Simons Foundation. The work of W.K. was partially supported by the Army Research Office under award W911NF1810306 and Air Force Office of Scientific Research under award number FA9550-23-1-0011. The authors thank Cameron Thieme for helpful discussions.

\bibliographystyle{plain}
\bibliography{GP_dynamics}

\end{document}